\documentclass[a4paper,10pt]{amsart}

\usepackage[colorlinks=printing]{hyperref}
\usepackage{amssymb,amsmath,esint,graphicx,subfigure}

\newtheorem{theorem}{Theorem}[section]
\newtheorem{lemma}[theorem]{Lemma}

\newtheorem{corollary}[theorem]{Corollary}
\newtheorem{definition}{Definition}[section]

\theoremstyle{remark} \newtheorem{remark}[theorem]{Remark} \theoremstyle{definition} 

\numberwithin{equation}{section}

\newcommand{\del}{\partial}

\newcommand{\ra}{\rightarrow}

\newcommand{\alp}{\alpha}
\newcommand{\eps}{\ensuremath{\varepsilon}}
\newcommand{\R}{\ensuremath{\mathbb{R}}}

\newcommand{\Z}{\ensuremath{\mathbb{Z}}}
\newcommand{\tg}{\ensuremath{\tilde g}}
\newcommand{\dx}{\ensuremath{\Delta x}}

\newcommand{\sgn}{{\rm sgn}\, }
\newcommand{\dif}{\mathrm{d}}

\begin{document}

\title[Degenerate fractional order convection-diffusion
equations]{Entropy solution theory for fractional degenerate
  convection-diffusion equations}

\author[{S.~Cifani}]{{Simone Cifani}}
\address[Simone Cifani]{\\ Department of Mathematics\\ Norwegian University of Science and Technology (NTNU)\\
 N-7491 Trondheim, Norway}
 \email[]{simone.cifani\@@math.ntnu.no}
\urladdr{http://www.math.ntnu.no/\~{}cifani/}

\author[{E.R.~Jakobsen}]{{Espen R. Jakobsen}}
\address[Espen R. Jakobsen]{\\ Department of Mathematics\\ Norwegian University of Science and Technology (NTNU)\\
 N-7491 Trondheim, Norway}
\email[]{erj\@@math.ntnu.no}
\urladdr{http://www.math.ntnu.no/\~{}erj/}



\thanks{This research was supported by the Research Council of Norway (NFR) through the project "Integro-PDEs:
Numerical methods, Analysis, and Applications to Finance".}

\begin{abstract}
  We study a class of degenerate convection diffusion equations with a
  fractional non-linear diffusion term. This class is a new, but natural,
  generalization of local degenerate convection diffusion equations, and include
  anomalous diffusion equations, fractional conservations laws,
  fractional Porous medium equations, and new fractional degenerate
  equations as special cases. We define
  weak  entropy solutions and prove
  well-posedness under weak regularity assumptions on the solutions,
  e.g.~uniqueness is obtained in the class of bounded integrable
  solutions. Then we introduce a new monotone conservative numerical
  scheme and prove convergence toward the entropy solution in the class
  of bounded integrable BV functions. The well-posedness results are
  then extended to non-local terms based on general L\'{e}vy
  operators, connections to some fully
  non-linear HJB equations are established, and finally, some numerical
  experiments are included to  give the reader an idea about the qualitative
  behavior of solutions of these new equations.
\end{abstract}

\maketitle

%\tableofcontents

\section{Introduction}

In this paper we study well-posedness and approximation of a Cauchy problem for the possibly degenerate non-linear non-local integral partial differential
equation
\begin{align}\label{1}
\begin{cases}
\partial_tu+\nabla\cdot f(u)=-(-\Delta)^{\lambda/2}A(u)&\text{in }Q_T=\mathbb{R}^d\times(0,T),\\
u(x,0)=u_{0}(x)&\text{in }\mathbb{R}^d,
\end{cases}
\end{align}
where $f=(f_1,\dots,f_d):\mathbb{R}\rightarrow\mathbb{R}$ and $A:\mathbb{R}\rightarrow\mathbb{R}$ are Lipschitz continuous with Lipschitz constants $L_f$ and
$L_A$, $A(\cdot)$ non-decreasing with $A(0)=0$, and the non-local operator $-(-\Delta)^{\lambda/2}$ (or $g[\cdot]$ in shorthand notation) is the fractional
Laplacian defined as
\begin{equation*}
-(-\Delta)^{\lambda/2}\phi(x)=c_\lambda\  P.V.\int_{|z|>0}\frac{\phi(x+z,t)-\phi(x,t)}{|z|^{d+\lambda}}\ \dif z
\end{equation*}
for some constants $c_\lambda>0$, $\lambda\in(0,2)$, and a sufficiently regular function $\phi$. Note that $A(\cdot)$ can be {\em strongly degenerate}, i.e.~it
may vanish on a set of positive measure.

Equation \eqref{1} is a fractional degenerate convection diffusion equation, and  this class of equations has received considerable interest recently thanks to
the wide variety of applications. They encompass various linear anomalous diffusion equations ($f\equiv0$ and $A(u)\equiv u$), scalar conservation laws
\cite{Da:Book,Holden/Risebro,Kruzkov,Oleinik,Se:Book} ($A\equiv0$), fractional (or \emph{fractal}) conservation laws \cite{Alibaud,Droniou/Imbert} ($A(u)\equiv
u$), and some (but not all!) fractional  Porous medium equations \cite{Depablo/Quiros/Rodriguez/Vasquez} ($f\equiv0$ and $A(u)=|u|u^m$, $m\geq 1$), but see
also \cite{Biler/Imbert/Karch,Caffarelli/Vazquez}. Equation \eqref{1} is an extension to the fractional diffusion setting of the degenerate
convection-diffusion equation \cite{Carrillo,Karlsen/Risebro}
\begin{equation}\label{conv_diff}
\begin{split}
\partial_tu+\nabla\cdot f(u)=\Delta A(u).
\end{split}
\end{equation}
When $A(\cdot)$ is strongly degenerate, equation \eqref{1} has never been analyzed before as far as we know.

The literature concerning the type of equations mentioned above is
immense. We will only give a partial and incomplete survey of some
parts we feel are more relevant for this paper. For a more complete
discussion and many more references, we refer the reader to the nice
papers \cite{Alibaud} and \cite{KU}. But before we continue, we would
like to mention actual and potential applications. A large variety of
phenomena in physics and finance are modeled by linear anomalous
diffusion equations, see
e.g. \cite{Woyczynski,Ap:Book,CT:Book}. Fractional conservation laws
are generalizations of convection-diffusion equations
(\eqref{conv_diff} with $A(u)\equiv u$), and appear in some physical
models for over-driven detonation in gases \cite{Clavin} and
semiconductor growth \cite{Woyczynski}, and in areas like dislocation
dynamics, hydrodynamics, and molecular biology,
cf. \cite{Alibaud,Alibaud/Droniou/Vovelle,Droniou}. Similar equations, but
with slightly different non local term, also appear in radiation hydrodynamics
\cite{RoYo}. Equations like \eqref{conv_diff} are
used to model a vast variety of phenomena, including porous media flow
\cite{Va:Book}, reservoir simulation \cite{EsKa}, sedimentation
processes \cite{BCBT}, and traffic flow \cite{Book:Wi}. Finally, we
mention \cite{JaWi} where degenerate elliptic-parabolic equations with
fractional time derivatives are considered.

In the non-linear and degenerate setting of \eqref{1}, we can not expect to have classical solutions and it is well-known that weak solutions are not unique in
general. In the setting of fractional conservation laws this is proved in e.g.~\cite{Alibaud/Andreianov,Alibaud/Droniou/Vovelle,Kiselev/Nazarov/Shterenberg}.
To get uniqueness we impose extra conditions, called entropy conditions. In this paper we will introduce a Kruzkov type entropy formulation for equation
\eqref{1}. This type of formulation was introduced by Kruzkov in \cite{Kruzkov}, and used along with a doubling of variables device, to obtain general
uniqueness results for scalar conservation laws. Much later, Carrillo in \cite{Carrillo} extended these results to cover second order equations like
\eqref{conv_diff}, see also \cite{Karlsen/Risebro} for more general results and a presentation and proof which is more like our own. More recently, Alibaud
\cite{Alibaud} extended the Kruzkov formulation and uniqueness result to the fractional setting. He obtained general results for fractional conservation laws.
In a new work by Karlsen and Ulusoy \cite{KU}, a unified formulation is given that essentially includes the results of Alibaud and Carrillo as special cases.
In \cite{Alibaud,KU} the fractional diffusion is always linear and non-degenerate.

The entropy formulation we use is an extension of the formulation of Alibaud, and it allows us to prove a general $L^1$-contraction and uniqueness result for
bounded integrable solutions of the initial value problem \eqref{1}. Our uniqueness proof relies on some new observations and estimates along with ideas from
\cite{Carrillo,Karlsen/Risebro}. From a technical point of view, our proof for $\lambda\in(0,2)$ is more related to the conservation law (or fractional
conservation law) proof than the  more technical proof of Carrillo for $\lambda=2$ (equation \eqref{conv_diff}). E.g.~we do not need a ``weak chain rule'' and
hence do not need to assume any extra a priori regularity on the term $A(u)$.

In practice to solve \eqref{1} we must resort to numerical computations. But since the equation is non-linear and degenerate, many numerical methods will fail
to converge or converge to false (non-entropy) solutions. The solution is to construct ``good'' numerical methods that insure convergence to entropy solutions.
In the conservation law community, it is well known that monotone, conservative, and consistent methods will do the job for you. There is a vast literature on
such methods, we refer the reader e.g.~to \cite{Holden/Risebro} and references therein. For non-linear fractional equations there exist very few methods and
results so far. Dedner and Rhode \cite{DeRo} introduced a convergent finite volume method for a non-local conservation laws from radiation hydrodynamics.
Droniou \cite{Droniou} was the first to define and prove convergence for approximations of fractional conservations laws. Karlsen and the authors then
introduced and proved convergence for Discontinuous Galerkin methods for fractional conservation laws and fractional convection-diffusion equations in
\cite{Cifani/Jakobsen/Karlsen,Cifani/Jakobsen/Karlsen2}. After that, the authors introduced a convergent spectral vanishing viscosity method for fractional
conservations laws in \cite{Cifani/Jakobsen}. Kuznetzov type error estimates were also obtained in \cite{Cifani/Jakobsen,Cifani/Jakobsen/Karlsen}. In this
paper, we discretize for the first time \eqref{1} in its general form. We introduce a new difference quadrature approximation that we prove converges to the
entropy solution. The convergence holds for bounded integrable BV solutions, and hence we also have existence of solutions
 in this class. Finally, existence of solutions in the wider class of bounded integrable function is obtained through approximation via bounded integrable BV solutions
(cf. Theorem \ref{th:existence}).

In many applications, especially in finance, the non-local term is not a fractional Laplacian, but rather a L\'{e}vy type operator $g_\mu$:
$$g_\mu[\phi](x)=\int_{|z|>0}\phi(x+z)-\phi(x)-z\cdot
\nabla\phi(x)1_{|z|<1}\ \dif\mu(z),$$ where the L\'{e}vy measure $\mu$ is a positive Radon measure satisfying
$$\int_{|z|>0}|z|^2\wedge 1 \ \mu(\dif z) <\infty.$$
These operators are the infinitesimal generators of pure jump L\'{e}vy processes. We refer to \cite{Ap:Book,CT:Book} for the theory and applications of such
processes and to \cite{KU} for a very relevant and nice discussion and many more references. The entropy solution theory related to such operators is very
similar to the one for fractional Laplacians, and the first well-posedness results were obtained in \cite{KU}. In this paper we extend the entropy theory for
\eqref{1} to this L\'{e}vy setting (cf. equation \eqref{1bis}). Our formulation is an extension of Alibaud's formulation and is different from the one given in
\cite{KU}. We also treat completely general L\'{e}vy measures, i.e.~our L\'{e}vy operators are slightly more general than the ones in \cite{KU}.

We also discuss the fact that \eqref{1} is related to fully non-linear HJB equations, see Section \ref{sec:HJB}. We first show an easy extension of results
from \cite{Li:Book}: In one space dimension the gradient of a viscosity solution of a fractional HJB equation is an entropy solution of a fractional
conservation law. Then we show a new correspondence for any space dimension: If $u$ is a viscosity solution of
\begin{equation*}
u_t - A(g_\mu[u])=0,
\end{equation*}
then $v=g_\mu[u]$ is the entropy solution of
\begin{equation*}
v_t -g_\mu[A(v)]=0.
\end{equation*}
The relevance of these results are discussed in Section
\ref{sec:HJB}. The final part of the paper is devoted to numerical
simulations to give the reader an idea about the qualitative behavior
of the solutions of these new equations.

Here is the content of the paper section by section. The entropy formulation is introduced and discussed in Section \ref{sec:def}. In Section \ref{sec:uniq},
we state and prove $L^1$-contraction and uniqueness for entropy solutions of \eqref{1}. The monotone conservative numerical method is then introduced and
analyzed in Section \ref{sec:num}. In Section \ref{sec:Levy}, we extend the well-posedness results proved for solutions of \eqref{1} to a wider class of
equations where the fractional Laplacian has been replaced by a general L\'{e}vy  operator. In  Section \ref{sec:HJB} we show how solutions of equations of the
type \eqref{1} are related to solutions of fully non-linear HJB equations, and in the last section, we provide several numerical simulations of problems of the
form \eqref{1}.

\section{Entropy formulation}
\label{sec:def} In this section we introduce an entropy formulation for the initial value problem \eqref{1} which generalizes Alibaud's formulation in
\cite{Alibaud}. To this end, let us split the non-local operator $g$ into two terms: for each $r>0$, we write $g[\varphi]=g_r[\varphi]+g^r[\varphi]$ where
\begin{align*}
g_r[\varphi](x)&=c_\lambda\ P.V.\int_{|z|<r}\frac{\varphi(x+z)-\varphi(x)}{|z|^{d+\lambda}}\ \dif z,\\
g^r[\varphi](x)&=c_\lambda\int_{|z|>r}\frac{\varphi(x+z)-\varphi(x)}{|z|^{d+\lambda}}\ \dif z.
\end{align*}
The Cauchy principal value is defined as
\begin{align*}
P.V.\int_{|z|>0}\varphi(z)\ \dif z=\lim_{b\ra0}\int_{b<|z|}\varphi(z)\ \dif z.
\end{align*}
Note that, by symmetry, $$P.V.\int_{|z|<r}\frac{z}{|z|^{d+\lambda}}\ \dif z=0$$ and hence
\begin{align}
\label{gdef2} g_r[\varphi](x)=c_\lambda\ P.V.\int_{|z|<r} \frac{\varphi(x+z)-\varphi(x)-z\cdot\nabla\varphi(x)}{|z|^{d+\lambda}}\ \dif z.
\end{align}
Whenever $\varphi$ is smooth enough, the principal value in
\eqref{gdef2} is well defined by the dominated convergence theorem
since
\begin{align*}
  |g_r[\varphi](x)|&\leq \left\{\begin{array}{ll}
\!\! c_\lambda\|D\varphi\|_{L^\infty(B(x,r))}\int_{|z|<r}\frac{|z|}{|z|^{d+\lambda}}\
  \dif z& \text{when }\lambda\in(0,1)\\[0.2cm]
\!\!  \frac{c_\lambda}{2}\|D^2\varphi\|_{L^\infty(B(x,r))}\int_{|z|<r}\frac{|z|^2}{|z|^{d+\lambda}}\
  \dif z& \text{when }\lambda\in[1,2)
\end{array}\right\}<\infty.
\end{align*}
The above integrals are finite because in
polar coordinates they are proportional to
$$\int_0^r
\frac{s^1}{s^{d+\lambda}}\ s^{d-1}\ \dif s\ \text{for }\lambda\in(0,1)\quad\text{and}\quad\int_0^r
\frac{s^2}{s^{d+\lambda}}\ s^{d-1}\ \dif s\ \text{for
}\lambda\in[1,2).$$
This estimate also shows that the integral in
\eqref{gdef2} exists and this leads to an alternative definition of
the operator $g_r$ avoiding the principal value (i.e.~\eqref{gdef2}
without $P.V.$). This second definition is used e.g.~in
\cite{Alibaud}.

Let us introduce the functions $\eta_{k}(u)=|u-k|$,
$\eta'_k(u)=\text{sgn}(u-k)$, and $q_{k}(u)=\eta'_k(u)(f(u)-f(k))$
where the sign function is defined as
$$\text{sgn}(s)=\begin{cases}1&\text{for }s>0\\
0&\text{for }s=0\\
-1&\text{for }s<0.\end{cases}$$
 The entropy
formulation we use is the following:

\begin{definition}\label{L1-entropy}
A function $u$ is an entropy solution of the initial value problem \eqref{1} provided that
\begin{itemize}
\item[\emph{i)}]$u\in L^{\infty}(Q_T)\cap C([0,T];L^{1}(\mathbb{R}^d))$;
\item[\emph{ii)}]for all $k\in\mathbb{R}$, all $r>0$, and all nonnegative test functions $\varphi\in C_{c}^{\infty}(Q_T)$,
\begin{equation*}
\begin{split}
\iint_{Q_T}\eta_k(u)\partial_t\varphi+q_k(u)\cdot\nabla\varphi+\eta_{A(k)}(A(u))\, g_r[\varphi]+\eta'_k(u)\, g^r[A(u)]\, \varphi\ \dif x\dif t\geq0;
\end{split}
\end{equation*}
\item[\emph{iii)}]$u(\cdot,0)=u_0(\cdot)$ a.e.
\end{itemize}
\end{definition}
\begin{remark}
By $C([0,T];L^1(\R^d))$ we mean the Banach space where the norm is given by $\|\phi\|_{C([0,T];L^1(\R^d))}=\max_{t\in[0,T]}\!\big\{\int_{\R^d}|\phi(x,t)|\ \dif
x\big\}.$
\end{remark}
\begin{remark}
\label{WDef} In view of  \emph{i)} and the properties of $f$ and $A$, $\eta_k(u),q_k(u),\eta_{A(k)}(A(u))\in L^\infty(Q_T)$ while $A(u)\in L^\infty(Q_T)\cap
C([0,T];L^1(\R^d))$. It immediately follows that the local terms in \emph{ii)} are well-defined. Since  $g_r[\phi]\in C_c^\infty(Q_T)$ for  $\phi\in
C_c^\infty(Q_T)$, also the $g_r$-term in  \emph{ii)} is well-defined. Finally we note that $g^r[\psi](x)$ is well-defined and belongs to $L^\infty(\R^d)$ for
$\psi\in L^\infty(\R^d)$, and to $L^1(\R^d)$ for $\psi\in L^1(\R^d)$ by Fubini (integrating first w.r.t.~$x$). It follows that  $g^r[A(u)]\in L^\infty(Q_T)\cap
C([0,T];L^1(\R^d))$, and hence that the $g^r$-term in \emph{ii)} is well-defined.
\end{remark}
\begin{remark}
Since $u\in C([0,T];L^1(\R^d))$ by \emph{i)}, part {\em iii)} implies that the initial condition is imposed in the strong $L^1$-sense:
$$\lim_{t\ra0}\|u(\cdot,t)-u_0\|_{L^1(\R^d)}=0.$$
A more traditional approach where initial values $u(\cdot,0)$ are
included in the entropy inequality \emph{ii)} would also work,
cf.~e.g.~\cite[Chapter 2]{Holden/Risebro}.
\end{remark}

Let us point out that, in the case $\lambda\in(0,1)$ and whenever the
entropy solutions are sought in the $BV$-class, Definition
\ref{L1-entropy} can be simplified to the following one:
\begin{definition}\label{BV-entropy}
A function $u$ is an entropy solution of the initial value problem \eqref{1} provided that
\begin{itemize}
\item[\emph{i)}]$u\in L^{\infty}(Q_T)\cap L^\infty(0,T;
  BV(\mathbb{R}^d))\cap C([0,T];L^1(\R^d))$;
\item[\emph{ii)}]for all $k\in\mathbb{R}$ and all nonnegative test functions $\varphi\in C_{c}^{\infty}(Q_T)$,
\begin{equation*}
\begin{split}
\iint_{Q_T}\eta_k(u)\partial_t\varphi+q_k(u)\cdot\nabla\varphi+\eta'_k(u)\, g[A(u)]\, \varphi\ \dif x\dif t\geq0;
\end{split}
\end{equation*}
\item[\emph{iii)}]$u(\cdot,0)=u_0(\cdot)$ a.e.
\end{itemize}
\end{definition}
Note that the non-local term $g[A(u)]$ in the integral in {\em ii)} is well defined as shown in the following lemma.
\begin{lemma}\label{interp}
If $\lambda\in(0,1)$, then there is a constant $C>0$ such that
\begin{equation*}
\begin{split}
\|g[A(u)]\|_{L^1(\R^d)}\leq c_\lambda CL_A\|u\|^{1-\lambda}_{L^1(\R^d)}|u|^\lambda_{BV(\R^d)}.
\end{split}
\end{equation*}
\end{lemma}

\begin{proof}We split the integral in two parts, use Fubini and the
  estimate
$$\int_{\R^d}|u(x+z)-u(x)|\ \dif x\leq \sqrt{d}|z||u|_{BV(\R^d)}$$
(cf.~Lemma \ref{lem:BV}), and change to polar coordinates ($z=ry$ for
$r\geq0$ and $|y|=1$) to find that:
\begin{align*}
&\int_{|z|<\epsilon}\int_{\R^d}\frac{|A(u(x+z))-A(u(x))|}{|z|^{d+\lambda}}\ \dif x\dif z\leq
L_A\sqrt d|u|_{BV(\R^d)}\int_{|z|<\eps}\frac{|z|}{|z|^{d+\lambda}}\ \dif z\\
&=L_A\sqrt d|u|_{BV(\R^d)}\int_{|y|=1}\dif S_y \int_0^\eps \frac{\dif r}{r^{\lambda}}=L_A\sqrt d|u|_{BV(\R^d)}\eps^{1-\lambda}\int_{|y|=1}\dif S_y\int_0^1
\frac{\dif r}{r^{\lambda}}
\end{align*}
and
\begin{align*}
\int_{|z|>\epsilon}\int_{\R^d}\frac{|A(u(x+z))-A(u(x))|}{|z|^{d+\lambda}}\ \dif x\dif z &\leq
\frac{2L_A}{\epsilon^{\lambda}}\|u\|_{L^{1}(\R^d)}\int_{|y|=1}\dif S_y\int_{1}^\infty\frac{\dif r}{r^{1+\lambda}}.
\end{align*}
To conclude, we choose $\epsilon=\|u\|_{L^{1}(\R^d)}|u|^{-1}_{BV(\R^d)}$.
\end{proof}

The following result shows how the two definitions of entropy
solutions are interrelated and how they relate to weak and classical solutions
of \eqref{1}.
\begin{theorem}$\quad$
\begin{itemize}
\item[\emph{i)}] Definition \ref{L1-entropy} and Definition \ref{BV-entropy} are equivalent whenever $\lambda\in(0,1)$
and $u\in L^{\infty}(Q_T)\cap L^\infty(0,T;
  BV(\mathbb{R}^d))\cap C([0,T];L^1(\R^d))$.\smallskip
\item[\emph{ii)}] Any entropy solution $u$ of \eqref{1} is a weak
  solution: for all $\varphi\in C_{c}^{\infty}(Q_T)$,
\begin{equation*}
\iint_{Q_T}u\partial_t\varphi+f(u)\cdot\nabla\varphi+A(u)\,g[\varphi]\ \dif x\dif t=0.
\end{equation*}

\item[\emph{iii)}] If $A\in C^2(\R)$, then any  classical solution $u\in L^{\infty}(Q_T)\cap  C([0,T];L^{1}(\mathbb{R}^d))$ of  \eqref{1}
is an entropy solution.
\end{itemize}
\end{theorem}

\begin{remark}
In $iii)$ we need additional regularity of $A$ to give a pointwise sense to the equation and hence also to define classical solutions. When $\lambda\in[1,2)$
it suffices to assume that $A\in C^2$, and when $\lambda\in(0,1)$ $A\in C^1$ is enough.
\end{remark}
\begin{proof}\*

\emph{i)} Repeated use of the dominated convergence theorem and Lemma \ref{interp} first shows that, when $r\rightarrow 0$,
$$g_r[\varphi]\rightarrow0\quad\text{and}\quad g^r[A(u)]\rightarrow
g[A(u)]\quad a.e., $$ and then combined with this convergence result and H\"older's inequality, that Definition \ref{L1-entropy} implies Definition
\ref{BV-entropy} when $u$ is BV. To go the other way, let us note that since $A(\cdot)$ is non-decreasing,
\begin{equation}\label{non_increasing}
\begin{split}
\sgn(u-k)(A(u)-A(k))=|A(u)-A(k)|.
\end{split}
\end{equation}
Thus, if we write
\begin{equation*}
\begin{split}
&g[A(u)]=g_\epsilon[A(u)]\\
&\quad+c_\lambda\int_{\epsilon<|z|<r}\frac{(A(u(x+z,t))-A(k))-(A(u(x,t))-A(k))}{|z|^{1+\lambda}}\ \dif z\\
&\quad+g^r[A(u)],
\end{split}
\end{equation*}
multiply each side by $\eta_k'(u)\varphi$ and integrate over $Q_T$, we end up with
\begin{equation*}
\begin{split}
&\iint_{Q_T}\eta_k'(u)\, g[A(u)]\, \varphi\ \dif x\dif t\leq\iint_{Q_T}\Bigg\{\eta_k'(u)\, g_\epsilon[A(u)]\, \varphi\\
&\quad+c_\lambda\varphi\int_{\epsilon<|z|<r}\frac{|A(u(x+z,t))-A(k)|-|A(u(x,t))-A(k)|}{|z|^{1+\lambda}}\ \dif z\\
&\quad+\eta_k'(u)\, g^r[A(u)]\, \varphi\Bigg\}\ \dif x\dif t.
\end{split}
\end{equation*}
We now use the change of variables $(z,x)\rightarrow(-z,x+z)$ to pass the test function $\varphi$ inside the integral $\epsilon<|z|<r$, and obtain
\begin{equation}\label{change_var}
\begin{split}
&\iint_{Q_T}\varphi(x,t)\int_{\epsilon<|z|<r}\frac{|A(u(x+z,t))-A(k)|-|A(u(x,t))-A(k)|}{|z|^{1+\lambda}}\ \dif z\dif x\dif t\\
&=\iint_{Q_T}|A(u(x,t))-A(k)|\int_{\epsilon<|z|<r}\frac{\varphi(x+z,t)-\varphi(x,t)}{|z|^{1+\lambda}}\ \dif z\dif x\dif t.
\end{split}
\end{equation}
The entropy inequality in Definition \ref{L1-entropy} is finally
recovered in the limit as $\epsilon\rightarrow0$.
\medskip

\emph{ii)} Using \eqref{non_increasing} and the change of variables $(z,x)\rightarrow(-z,x+z)$,
\begin{equation*}
\begin{split}
&\iint_{Q_T}\eta'_k(u(x,t))\, g^r[A(u(x,t))]\, \varphi(x,t)\ \dif x\dif t\\
&\leq c_\lambda\iint_{Q_T}\varphi(x,t)\int_{|z|>r}\frac{|A(u(x+z,t))-A(k)|-|A(u(x,t))-A(k)|}{|z|^{1+\lambda}}\ \dif z\dif x\dif t\\
&=\iint_{Q_T}|A(u(x,t))-A(k)|\, g^r[\varphi(x,t)]\, \dif x\dif t.
\end{split}
\end{equation*}
Thus, since $g=g_r+g^r$, we have produced the inequality
\begin{equation*}
\begin{split}
\iint_{Q_T}\eta_k(u)\partial_t\varphi+q_k(u)\cdot\nabla\varphi+\eta_{A(k)}(A(u))\,g[\varphi]\ \dif x\dif t\geq0.
\end{split}
\end{equation*}
By this inequality and the definitions of $\eta$ and $q$, if $\pm
k\geq\|u\|_{L^\infty(\R)}$, then
\begin{align*}
\mp \iint_{Q_T}(u-k)\del_t\phi+(f(u)-f(k))\cdot\nabla\phi+(A(u)-A(k))g[\phi]\ \dif x\dif t\geq0.
\end{align*}
By the Divergence theorem and a computation like in \eqref{change_var}, all
the $k$-terms are zero and hence $u$ is a weak solution as defined in
\emph{ii)}.
\medskip

\emph{iii)} Since $u$ solves equation \eqref{1} point-wise, for each $(x,t)\in Q_T$ and all $k\in\mathbb{R}$, we can write
\begin{equation*}
\begin{split}
&\partial_t(u-k)+\nabla\cdot(f(u)-f(k))=g_\epsilon[A(u)]\\
&+c_\lambda \int_{\epsilon<|z|<r}\frac{ (A(u(x+z,t))-A(k))-(A(u(x,t))-A(k)) }{|z|^{d+\lambda}}\ \dif z\\
&+g^r[A(u)].
\end{split}
\end{equation*}
If we multiply both sides of this equation by $\eta'_k(u)$ and use
\eqref{non_increasing}, we obtain
\begin{equation*}
\begin{split}
&\eta'_k(u)\,\partial_t(u-k)+\eta'_k(u)\,\nabla\cdot(f(u)-f(k))\leq\eta'_k(u)\,g_\epsilon[A(u)]\\
&+c_\lambda \int_{\epsilon<|z|<r}\frac{ |A(u(x+z,t))-A(k)|-|A(u(x,t))-A(k)|}{|z|^{d+\lambda}}\ \dif z\\
&+\eta'_k(u)\,g^r[A(u)].
\end{split}
\end{equation*}
Let us now multiply both sides of this inequality by a nonnegative test function $\varphi$, and integrate over $Q_T$ to obtain
\begin{equation*}
\begin{split}
&-\iint_{Q_T}\eta_k(u)\,\partial_t\varphi+q_k(u)\cdot\nabla\varphi\ \dif x\dif t\\
&\leq\iint_{Q_T}\Bigg\{\eta'_k(u)\, g_\epsilon[A(u(x,t))]\ \varphi\\
&\qquad+c_\lambda \varphi\int_{\epsilon<|z|<r}\frac{ |A(u(x+z,t))-A(k)|-|A(u(x,t))-A(k)| }{|z|^{d+\lambda}}\ \dif z\\
&\qquad+\eta'_k(u)\, g^r[A(u(x,t))]\, \varphi\Bigg\}\ \dif x\dif t.
\end{split}
\end{equation*}
Thanks to \eqref{change_var}, we can pass the test function $\varphi$ inside the integral $\epsilon<|z|<r$, and so recover the entropy inequality in Definition
\ref{L1-entropy} in the limit as $\epsilon\rightarrow0$.
\end{proof}

%%%%%%%%%%%%%%%%%%%%%%%%%%%%%%%%%%%%%%%%%%%%%%%%%%%%%%%%%%%%%%%%%%%%%%%%%%%%%%%%%%%%%%%%%%%%%%%%%%%%%%%%%%%%%%%%%%%%%%%%%%%%%%%%%%%%%%%%%%%%%%%%%%%%%%%%%%%%%%%%%%%%%%%%%%%
%%%%%%%%%%%%%%%%%%%%%%%%%%%%%%%%%%%%%%%%%%%%%%%%%%%%%%%%%%%%%%%%%%%%%%%%%%%%%%%%%%%%%%%%%%%%%%%%%%%%%%%%%%%%%%%%%%%%%%%%%%%%%%%%%%%%%%%%%%%%%%%%%%%%%%%%%%%%%%%%%%%%%%%%%%%

\section{$L^1$-contraction and Uniqueness}
\label{sec:uniq} We now establish $L^1$-contraction and uniqueness for entropy solutions of the initial value problem \eqref{1} using the Kru\v{z}kov's
doubling of variables device \cite{Kruzkov}. This technique has already been extended to fractional conservation laws (i.e., $A(u)=u$) by Alibaud
\cite{Alibaud}. The first part of our proof builds on the ideas developed by Alibaud (and Kru\v{z}kov!), but in the rest of the proof different ideas have to
be used in our non-linear and possibly degenerate setting.

\begin{theorem}\label{L1contraction}
  Let $u$ and $v$ be two entropy solutions of the initial value
  problem \eqref{1} with initial data $u_0$ and $v_0$. Then, for
   all $t\in(0,T)$,
\begin{equation*}
\|u(\cdot,t)-v(\cdot,t)\|_{L^1(\mathbb{R}^d)}\leq\|u_0-v_0\|_{L^1(\mathbb{R}^d)}.
\end{equation*}
\end{theorem}
Uniqueness for entropy solutions of \eqref{1} immediately follows from the above $L^1$-contraction: if $u_0=v_0$, then $u=v$ a.e.~on $Q_T$.
\begin{corollary}\emph{(Uniqueness)}
There is at most one entropy solution of \eqref{1}.
\end{corollary}

\begin{proof}[Proof of Theorem \ref{L1contraction}]\*

  1)$\quad$We take $u=u(x,t)$ and $v=v(y,s)$, let $\psi=\psi(x,y,t,s)$
  be a nonnegative test function, and denote by $\eta(u,k)$, $q(u,k)$,
  $\eta'(u,k)$ the quantities $\eta_k(u)$, $q_k(u)$,
  $\eta'_k(u)$. After integrating the entropy inequality for
  $u=u(x,t)$ with $k=v(y,s)$ over $(y,s)\in Q_T$, we find that
\begin{equation}\label{a1}
\begin{split}
\iint_{Q_T}\iint_{Q_T}&\eta(u(x,t),v(y,s))\, \partial_t\psi(x,y,t,s)\\
&+q(u(x,t),v(y,s))\cdot\nabla_x\psi(x,y,t,s)\\
&+\eta(A(u(x,t)),A(v(y,s)))\ g_r[\psi(\cdot,y,t,s)](x)\\
&+\eta'(u(x,t),v(y,s))\ g^r[A(u(\cdot,t))](x)\ \psi(x,y,t,s)\ \dif x\dif t\dif y\dif s\geq0.
\end{split}
\end{equation}
Similarly, since $\eta(u,k)=\eta(k,u)$, $q(u,k)=q(k,u)$, and
$\eta'(u,k)=-\eta'(k,u)$, integrating the entropy inequality for
$v=v(y,s)$ with $k=u(x,t)$ leads to
\begin{equation}\label{a2}
\begin{split}
\iint_{Q_T}\iint_{Q_T}&\eta(u(x,t),v(y,s))\, \partial_s\psi(x,y,t,s)\\
&+q(u(x,t),v(y,s))\cdot\nabla_y\psi(x,y,t,s)\\
&+\eta(A(u(x,t)),A(v(y,s)))\ g_r[\psi(x,\cdot,t,s)](y)\\
&-\eta'(u(x,t),v(y,s))\ g^r[A(v(\cdot,s))](y)\ \psi(x,y,t,s)\ \dif y\dif s\dif x\dif t\geq0.
\end{split}
\end{equation}
Let us now introduce the operator
\begin{equation*}
\begin{split}
\tg^r[\varphi(\cdot,\cdot)](x,y)=\int_{|z|>r}\frac{\varphi(x+z,y+z)-\varphi(x,y)}{|z|^{d+\lambda}}\ \dif z.
\end{split}
\end{equation*}
Since all the terms in \eqref{a1}--\eqref{a2} are integrable, we are are free to change the order of integration, and hence add up inequalities
\eqref{a1}--\eqref{a2} to find that (from now on $\dif w=\dif x\,\dif t\,\dif y\,\dif s$)
\begin{equation}\label{julie6}
\begin{split}
\iint_{Q_T}\iint_{Q_T}&\eta(u(x,t),v(y,s))\, (\partial_t+\partial_s)\psi(x,y,t,s)\\
&+q(u(x,t),v(y,s))\cdot(\nabla_x+\nabla_y)\psi(x,y,t,s)\\
&+\eta(A(u(x,t)),A(v(y,s)))\ g_r[\psi(\cdot,y,t,s)](x)\\
&+\eta(A(u(x,t)),A(v(y,s)))\ g_r[\psi(x,\cdot,t,s)](y)\\
&+\eta'(u(x,t),v(y,s))\ \tilde g^r[A(u(\cdot,t))-A(v(\cdot,s))](x,y)\ \psi(x,y,t,s)\ \dif w\geq0.
\end{split}
\end{equation}
In the following we will manipulate the operator $\tg^r$, while the
operators $g_r$ will simply be carried along to finally  vanish
in the limit as $r\rightarrow0$.

Let us use \eqref{non_increasing} to obtain the (Kato type of) inequality
\begin{equation*}
\begin{split}
&\eta'(u(x,t),v(y,s))\Big[\Big(A(u(x+z,t))-A(v(y+z,s))\Big)-\Big(A(u(x,t))-A(v(y,s))\Big)\Big]\\
&\leq |A(u(x+z,t))-A(v(y+z,s))|-|A(u(x,t))-A(v(y,s))|,
\end{split}
\end{equation*}
which implies that
\begin{equation}\label{j1}
\begin{split}
\eta'(u(x,t),v(y,s))\ \tilde g^r[A(u(\cdot,t))-A(v(\cdot,s))](x,y)\leq \tilde g^r\Big[|A(u(\cdot,t))-A(v(\cdot,s))|\Big](x,y).
\end{split}
\end{equation}
Furthermore, we use Fubini's Theorem and the change of variables
$(z,x,y)\rightarrow(-z,x+z,y+z)$  to see that
\begin{equation}\label{j2}
\begin{split}
&\iint_{Q_T}\iint_{Q_T}\psi(x,y,t,s)\ \tilde g^r\Big[|A(u(\cdot,t))-A(v(\cdot,s))|\Big](x,y)\ \dif w\\
&=\iint_{Q_T}\iint_{Q_T}|A(u(x,t))-A(u(y,s))|\ \tilde g^r[\psi(\cdot,\cdot,t,s)](x,y)\ \dif w.
\end{split}
\end{equation}
To sum up, when used in \eqref{julie6}, \eqref{j1}--\eqref{j2} produce the inequality
\begin{equation}\label{a3}
\begin{split}
\iint_{Q_T}\iint_{Q_T}&\eta(u(x,t),v(y,s))\, (\partial_t+\partial_s)\psi(x,y,t,s)\\
&+q(u(x,t),v(y,s))\cdot(\nabla_x+\nabla_y)\psi(x,y,t,s)\\
&+\eta(A(u(x,t)),A(v(y,s)))\ g_r[\psi(\cdot,y,t,s)](x)\\
&+\eta(A(u(x,t)),A(v(y,s)))\ g_r[\psi(x,\cdot,t,s)](y)\\
&+\eta(A(u(x,t)),A(v(y,s)))\ \tilde g^r[\psi(\cdot,\cdot,t,s)](x,y)\ \dif w\geq0.
\end{split}
\end{equation}
Thanks to the regularity of the test function $\psi$, we can now take
the limit as $r\rightarrow 0$ in \eqref{a3}, and end up with
\begin{equation}\label{ineq}
\begin{split}
\iint_{Q_T}\iint_{Q_T}&\eta(u(x,t),v(y,s))\, (\partial_t+\partial_s)\psi(x,y,t,s)\\
&+q(u(x,t),v(y,s))\cdot(\nabla_x+\nabla_y)\psi(x,y,t,s)\\
&+\eta(A(u(x,t)),A(v(y,s)))\ \tilde g[\psi(\cdot,\cdot,t,s)](x,y)\ \dif w\geq0,
\end{split}
\end{equation}
where
\begin{equation*}
\begin{split}
\tg[\varphi(\cdot,\cdot)](x,y)=P.V.\int_{|z|>0}\frac{\varphi(x+z,y+z)-\varphi(x,y)}{|z|^{d+\lambda}}\ \dif z.
\end{split}
\end{equation*}
Inequality \eqref{ineq} concludes the first part of the proof.
\medskip

2)$\quad$We now specify the test function $\psi$ in order to
derive the $L^1$-contraction from inequality \eqref{ineq}:
\begin{equation*}
\begin{split}
\psi(x,t,y,s)=\hat\omega_\rho\left(\frac{x-y}{2}\right)\omega_{\rho}\left(\frac{t-s}{2}\right)\phi\left(\frac{x+y}{2},\frac{t+s}{2}\right),
\end{split}
\end{equation*}
for $\rho>0$ and some $\phi\in C_c^\infty(Q_T)$ to be chosen later. Here $\hat\omega_\rho(x)=\omega_\rho(x_1)\cdots \omega_\rho(x_d)$ and
$\omega_\rho(s)=\frac{1}{\rho}\omega(\frac{s}{\rho})$ for a nonnegative $\omega\in C_c^\infty(\mathbb{R})$ satisfying
\begin{equation*}
\begin{split}
\omega(-s)=\omega(s),\quad \omega(s)=0 \text{ for all $|s|\geq 1$,$\quad$and}\quad \int_{\mathbb{R}}\omega(s)\ \dif s=1.
\end{split}
\end{equation*}
The reader can easily check that
\begin{align*}
(\partial_t+\partial_s)\psi(x,y,t,s)
&=\hat\omega_\rho\Big(\frac{x-y}{2}\Big)\omega_{\rho}\Big(\frac{t-s}{2}\Big)(\partial_t+\partial_s)\phi\Big(\frac{x+y}{2},\frac{t+s}{2}\Big),\\
(\nabla_x+\nabla_y)\psi(x,y,t,s)
&=\hat\omega_\rho\Big(\frac{x-y}{2}\Big)\omega_{\rho}\Big(\frac{t-s}{2}\Big)(\nabla_x+\nabla_y)\phi\Big(\frac{x+y}{2},\frac{t+s}{2}\Big),\\
\tg[\psi(\cdot,\cdot,t,s)](x,y)&=\hat\omega_\rho\Big(\frac{x-y}{2}\Big)\omega_{\rho}\Big(\frac{t-s}{2}\Big)
g\Big[\phi\Big(\cdot,\frac{t+s}{2}\Big)\Big]\left(\frac{x+y}{2}\right).
\end{align*}
Note that with this choice of test function $\psi$, expressions involving $\tg$ naturally transform into expressions involving $g$.

We now show that, in the limit $\rho\rightarrow0$, inequality \eqref{ineq} reduces to
\begin{equation}\label{julie10}
\begin{split}
\iint_{Q_T}&\eta(u(x,t),v(x,t))\partial_t\phi(x,t)\\
&+q(u(x,t),v(x,t))\cdot\nabla\phi(x,t)\\
&+\eta(A(u(x,t)),A(v(x,t))\ g[\phi(\cdot,t)](x)\ \dif x\dif t\geq0.
\end{split}
\end{equation}
Loosely speaking the reason for this is that the function $\omega_\delta$
 converges to the $\delta$-measure. A proof concerning
the local terms can be
found in e.g.~\cite{Karlsen/Risebro}. It remains to prove that
\begin{align*}
M:=&\Bigg|\iint_{Q_T}\iint_{Q_T}|A(u(x,t))-A(v(y,s))|\\
&\qquad\qquad\qquad\hat\omega_\rho\left(\frac{x-y}{2}\right)\omega_{\rho}\left(\frac{t-s}{2}\right)g\Big[\phi\Big(\cdot,\frac{t+s}{2}\Big)\Big]\left(\frac{x+y}{2}\right)\
\dif w\\
&\qquad\qquad-\iint_{Q_T}|A(u(x,t))-A(v(x,t))|\ g[\phi(\cdot,t)](x)\ \dif x\dif t\Bigg|\stackrel{\rho\rightarrow0}{\longrightarrow}0.
\end{align*}
To see this, we add and subtract
\begin{equation*}
\begin{split}
&\iint_{Q_T}\iint_{Q_T}|A(u(x,t))-A(v(x,t))|\\
&\qquad\qquad\qquad\hat\omega_\rho\left(\frac{x-y}{2}\right)\omega_{\rho}\left(\frac{t-s}{2}\right)
g\Big[\phi\Big(\cdot,\frac{t+s}{2}\Big)\Big]\left(\frac{x+y}{2}\right)\ \dif w,
\end{split}
\end{equation*}
use the fact that $\iint_{Q_T}\hat\omega_\rho\left(\frac{x-y}{2}\right)\omega_{\rho}\left(\frac{t-s}{2}\right)\dif y \dif s=1$ for any fixed $t\in(0,T)$ for
$\rho$ small enough, and that $\phi$ has compact support in $(0,T)$ to find that
\begin{align*}
M\leq&\iint_{Q_T}\iint_{Q_T}\Big||A(u(x,t))-A(v(y,s))|-|A(u(x,t))-A(v(x,t))|\Big|\\
&\qquad\qquad\qquad\qquad
\hat\omega_\rho\left(\frac{x-y}{2}\right)\omega_{\rho}\left(\frac{t-s}{2}\right)g\Big[\phi\Big(\cdot,\frac{t+s}{2}\Big)\Big]\left(\frac{x+y}{2}\right)\
\dif w\\[0.2cm]
&+\iint_{Q_T}\iint_{Q_T}\Bigg|g\Big[\phi\Big(\cdot,\frac{t+s}{2}\Big)\Big]\left(\frac{x+y}{2}\right)-g[\phi(\cdot,t)](x)\Bigg|\\
&\qquad\qquad\qquad\qquad\hat\omega_\rho\left(\frac{x-y}{2}\right)\omega_{\rho}\left(\frac{t-s}{2}\right)|A(u(x,t))-A(v(x,t))|\ \dif w.
\end{align*}
Let $M_1$ and $M_2$ denote the two integrals on the right hand side of
the expression above. By the inequality $||a-c|-|b-c||\leq|a-b|$ we see that
\begin{align*}
M_1\leq K_\phi\iint_{Q_T}\iint_{Q_T}|A(v(x,t))-A(v(y,s))|\ \hat\omega_\rho\left(\frac{x-y}{2}\right)\omega_{\rho}\left(\frac{t-s}{2}\right)\dif w,
\end{align*}
since, for all $(x,t),(y,s)\in Q_T\times Q_T$,
\begin{equation}\label{bound}
\begin{split}
&\Big|g\Big[\phi\Big(\cdot,\frac{t+s}{2}\Big)\Big]\left(\frac{x+y}{2}\right)\Big|\\
&\leq K_\phi:=\frac{c_\lambda}{2}\|D^2\phi\|_{L^\infty(\mathbb{R}^d)}\int_{|z|<1}\frac{|z|^2}{|z|^{d+\lambda}}\ \dif
z+2c_\lambda\|\phi\|_{L^\infty(\mathbb{R})}\int_{|z|>1}\frac{\dif z}{|z|^{d+\lambda}}.
\end{split}
\end{equation}
Note that both integrals in \eqref{bound} are finite (use polar coordinates to see this). Using the change of variables $x-y=h$ and $t-s=\tau$, we obtain
\begin{align*}
&M_1\\
&\leq K_\phi\iint_{Q_T}\iint_{Q_T}|A(v(x,t))-A(v(x+h,t+\tau))|\ \hat\omega_\rho\left(\frac{h}{2}\right)\omega_{\rho}\left(\frac{\tau}{2}\right)\dif x\dif
t\dif h\dif \tau\\
&\leq K_\phi\iint_{Q_T}\hat\omega_\rho\left(\frac{h}{2}\right)\omega_{\rho}\left(\frac{\tau}{2}\right)\left(\iint_{Q_T}|A(v(x,t))-A(v(x+h,t+\tau))|\ \dif x\dif
t\right)\dif h\dif \tau\\
&\leq K_\phi \sup_{|h|,|\tau|\leq \rho}\left(\iint_{Q_T}|A(v(x,t))-A(v(x+h,t+\tau))|\ \dif x\dif t\right)\stackrel{\rho\ra\infty}{\longrightarrow} 0
\end{align*}
by continuity of translations in $L^1$. We refer to Lemma 2.7.2 in
\cite{Se:Book} for a similar proof. A similar
argument using the fact that $g[\phi]\in C([0,T];L^1(\R^d))$
(cf. Remark \ref{WDef}) shows that $M_2\ra0$ as $\rho\ra0$, and we can
  therefore conclude that $M\leq M_1+M_2\ra0$ as $\rho\ra0$. The proof
  of \eqref{julie10} is now complete.
\medskip

3)$\quad$We now show that inequality \eqref{julie10} can be reduced to
\begin{align}\label{julie20}
\iint_{Q_T}|u(x,t)-v(x,t)|\,\chi'(t)\ \dif x\dif t\geq0,
\end{align}
if we take $\phi=\varphi_r(x)\chi(t)$ and send $r\ra\infty$ for $r>1$, $\chi\in C^\infty_c(0,T)$ (with derivative $\chi'$) to be specified later, and
\begin{align*}
\varphi_r(x)&=\int_{\mathbb{R}^d}\hat\omega(x-y)\mathbf{1}_{|y|<r}\ \dif y.
\end{align*}
All derivatives of $\varphi_r$ are bounded uniformly in $r$ and vanish for all $||x|-r|>1$. Concerning the flux-term in \eqref{julie10}, we find that
\begin{align*}
&\iint_{Q_T}\text{sgn}(u(x,t)-v(x,t))(f(u(x,t))-f(v(x,t)))\cdot\nabla\phi(x,t)\ \dif x\dif t\\
&\leq L_f\|\chi\|_{L^\infty}\iint_{Q_T}\Big(|u(x,t)|+|v(x,t)|\Big)\mathbf{1}_{||x|-r|<1}\ \dif x\dif t\stackrel{r\ra\infty}{\longrightarrow}0
\end{align*}
by the dominated convergence theorem since $u$ and $v$ belong to $L^1$ and $\mathbf{1}_{||x|-r|<1}\rightarrow 0$ as $r\ra\infty$ for all $x\in\mathbb{R}^d$.
The term in \eqref{julie10} containing the non-local operator also tends to zero as $r\ra\infty$. To see this note that $|g[\varphi_r](x)|$ is uniformly
bounded in $r$, cf.~\eqref{bound}, so by integrability of
$u$ and $v$ and H\"older's inequality,
 \begin{align*}
 &\iint_{Q_T}|A(u(x,t))-A(v(x,t))|\ |g[\varphi_r](x)|\ \dif x\dif t\\
 &\leq L_A\Big(\|u\|_{L^1(Q_T)}+\|v\|_{L^1(Q_T)}\Big)\sup_{r>1}\|g[\varphi_r]\|_{L^\infty(Q_T)}<\infty.
 \end{align*}
Hence we find that the integrand is bounded by an $L^1$-function
uniformly for $r>1$:
$$|A(u(x,t))-A(v(x,t))| |g[\varphi_r](x)|\leq L_A|(u(x,t)-v(x,t)| \sup_{r>1}\|g[\varphi_r]\|_{L^\infty(Q_T)}.$$
Then for any $x,z\in\R^d$ fixed
and $r>|x|+1$, $\varphi_r(x)=1$ and
$$|\varphi_r(x+z)-\varphi_r(x)|\leq
|\mathbf1_{|x+z|<r-1}-1|\leq\mathbf1_{|z|>r-1-|x|}.$$
With this in mind we find that
\begin{equation*}
|g[\varphi_r](x)|\leq
\int_{|z|>0}\frac{\mathbf1_{|z|>r-1-|x|}}{|z|^{d+\lambda}}\ \dif z\stackrel{r\ra\infty}{\longrightarrow} 0,
\end{equation*}
and hence  we can conclude by the dominated convergence theorem that
\begin{equation*}
\lim_{r\rightarrow\infty}\iint_{Q_T}|A(u(x,t))-A(v(x,t))|\ |g[\varphi_r](x)|\ \dif x\dif t=0.
\end{equation*}
\medskip

4)$\quad$To conclude the proof, we now take $\chi=\chi_\mu$ for
\begin{align*}
\chi_\mu(t)&=\int_{-\infty}^{t}(\omega_\mu(\tau-t_1)-\omega_\mu(\tau-t_2))\ \dif\tau,
\end{align*}
where $r>1$ and $0<t_1<t_2<T$. Loosely speaking, the function $\chi_\mu$ is a smooth approximation of the indicator function $\mathbf{1}_{(t_1,t_2)}$ which is
zero near $t=0$ and $t=T$ when $\mu>0$ is small enough. Since $\chi_\mu'(t)=\omega_\mu(t-t_1)-\omega_\mu(t-t_2)$, inequality \eqref{julie20} reduces to
\begin{equation*}
\iint_{Q_T}|u(x,t)-v(x,t)|\,\omega_\mu(t-t_2)\ \dif x\dif t\leq \iint_{Q_T}|u(x,t)-v(x,t)|\,\omega_\mu(t-t_1)\ \dif x\dif t.
\end{equation*}
By taking $\mu$ small enough and using Fubini's theorem, we can
rewrite this inequality as
\begin{equation}\label{mmm}
\Phi*\omega_\mu(t_2)\leq \Phi*\omega_\mu(t_1)\qquad \text{for}\qquad \Phi(t)=\int_{\R^d}|u(x,t)-v(x,t)|\ \dif x,
\end{equation}
where $\phi_1*\phi_2(t)=\int_{\R}\phi_1(s)\,\phi_2(t-s)\,\dif
s$. Since $u,v\in C([0,T];L^1(\R^d))$, we see that $\Phi\in C([0,T])$,
and hence by standard properties of convolutions,
\begin{equation*}
\Phi*\omega_\mu(t)\ra\Phi(t) \quad\text{as}\quad \mu\ra0.
\end{equation*}
for all $t\in(0,T)$. Hence we can send $\mu\ra0$ in \eqref{mmm} to obtain
\begin{equation*}
\|(u-v)(\cdot,t_2)\|_{L^1(\mathbb{R}^d)}\leq\|(u-v)(\cdot,t_1)\|_{L^1(\mathbb{R}^d)} .
\end{equation*}
Finally, the theorem follows from renaming $t_2$ and sending
$t_1\rightarrow0$ using $iii)$ and $C([0,T];L^1(\R^d))$ regularity of
$u$ and $v$.
\end{proof}

%%%%%%%%%%%%%%%%%%%%%%%%%%%%%%%%%%%%%%%%%%%%%%%%%%%%%%%%%%%%%%%%%%%%%%%%%%%%%%%%%%%%%%%%%%%%%%%%%%%%%%%%%%%%%%%%%%%%%%%%%%%%%%%%%%%%%%%%%%%%%%%%%%%%%%%%%%%%%%%%%%%%%%%%%%%
%%%%%%%%%%%%%%%%%%%%%%%%%%%%%%%%%%%%%%%%%%%%%%%%%%%%%%%%%%%%%%%%%%%%%%%%%%%%%%%%%%%%%%%%%%%%%%%%%%%%%%%%%%%%%%%%%%%%%%%%%%%%%%%%%%%%%%%%%%%%%%%%%%%%%%%%%%%%%%%%%%%%%%%%%%%

\section{A convergent numerical method}
\label{sec:num} In this section we introduce a numerical method for the initial value problem \eqref{1} which is monotone and conservative. Then we prove that
the limit of any convergent sequence of solutions of the method (as $\Delta x\ra0$) is an entropy solution of \eqref{1}. Finally we prove that any sequence of
solutions of the method  is relatively compact whenever the initial datum is a bounded integrable function of bounded variation, and hence we establish the
existence of an entropy solution of \eqref{1} in this case. Some numerical simulations based on this method are presented in the last section.

\subsection{Definition and properties of the numerical method}
For simplicity we only consider uniform space/time grids and we start
by the one dimensional case. The spatial grid then consists of the
points $x_i=i\Delta x$ for $i\in\Z$ and the
temporal grid of $t_n=n\Delta t$ for $n=0,\ldots,N$ and $N\Delta t=T$.
The explicit numerical method we consider then takes the form
\begin{equation*}
\left\{
\begin{split}
U_{i}^{n+1}&=U_{i}^{n}-\Delta t\,D^-F(U_{i}^{n},U_{i+1}^{n})+\Delta t\sum_{j\neq 0}G_j(A(U^{n}_{i+j})-A(U^{n}_{i})),\\
U_{i}^{0}&=\frac{1}{\Delta x}\int_{x_i+\Delta x[0,1)}u_{0}(x)\ \dif x,
\end{split}
\right.
\end{equation*}
where $D^-U_i=\frac{1}{\Delta x}(U_i-U_{i-1})$, $F:\mathbb{R}^2\rightarrow\mathbb{R}$ is a numerical flux satisfying
\begin{itemize}
\item[\emph{a)}] $F$ is Lipschitz continuous with Lipschitz constant $L_F$,
\item[\emph{b)}] $F$ is consistent,  $F(u,u)=f(u)$ for all
  $u\in\mathbb{R}$,
\item[\emph{c)}] $F(u_1,u_2)$ is non-decreasing w.r.t. $u_1$ and
  non-increasing w.r.t. $u_2$,
\end{itemize}
and $G_i$ is defined by
\begin{equation*}
\begin{split}
G_i&=c_\lambda\int_{x_i+\frac{\Delta x}2[-1,1)}\frac{\dif z}{|z|^{1+\lambda}}\quad\text{for $i\neq0$.}
\end{split}
\end{equation*}
In the multi dimensional case the spatial grid is $\Delta x\,\Z^d$
($\Delta x>0)$ with points
$$x_\alp=\Delta x\, \alp \qquad\text{where $\alp=(\alp_1,\dots,\alp_d)\in \Z^d$.}$$
Let $e_l$ be the $d$-vector with $l$-component 1 and the other components 0
and define the two box domains
$$R={\Delta x} [0,1)^d\quad \text{and}\quad R_0=\frac{\Delta x}2
[-1,1)^d,$$ noting that $\cup_\alp(x_\alp+R)=\cup_\alp(x_\alp+R_0)=\R^d$. The explicit numerical method we consider now takes the form
\begin{equation}\label{scheme}
\left\{
\begin{split}
U_{\alp}^{n+1}&=U_{\alp}^{n}-\Delta t\sum_{l=1}^dD^-_lF_l(U_{\alp}^{n},U_{\alp+e_l}^{n})+\Delta t\sum_{\beta\neq 0}G_\beta(A(U^{n}_{\alp+\beta})-A(U^{n}_{\alp})),\\
U_{\alp}^{0}&=\frac{1}{\Delta x^d}\int_{x_\alp+R}u_{0}(x)\ \dif x,
\end{split}
\right.
\end{equation}
where $D_l^-U_\alp=\frac{1}{\Delta x}(U_\alp-U_{\alp-e_l})$, $F_l:\mathbb{R}^2\rightarrow\mathbb{R}$ is a numerical flux satisfying \emph{a)} -- \emph{c)}
above with $f_l$ replacing $f$, and $G_\alp$ is defined by
\begin{equation*}
\begin{split}
G_\alp&=c_\lambda\int_{x_\alp+R_0}\frac{\dif z}{|z|^{d+\lambda}}\quad\text{for $\alp\neq0$.}
\end{split}
\end{equation*}
Note that $G_\alp$ is positive and finite since $0\not\in x_\alp+R_0$ unless $\alp=0$.

\begin{remark}
An admissible numerical flux $F_l$ is e.g.~the Lax-Friedrichs flux,
$$F_l(U_{\alp}^{n},U_{\alp+e_l}^{n})=\frac{1}{2}\left(f(U_{\alp}^{n})+f(U_{\alp+e_l})-\frac{\Delta x}{\Delta t}(U_{\alp+e_l}-U_{\alp}^{n})\right).$$
We refer the reader to \cite{Eymard} or \cite[Chapter 3]{Holden/Risebro} for a detailed presentation of more numerical fluxes which fulfill assumptions
\emph{a) -- c)}.
\end{remark}

Let us introduce the piecewise constant space/time interpolation
\begin{equation*}
\bar{u}(x,t)=U_{\alp}^{n}\quad\text{for all $(x,t)\in(x_\alp+R)\times[t_n,t_{n+1})$.}
\end{equation*}
In the following we often need the relation
\begin{equation}\label{scheme_fund_rel}
\begin{split}
\sum_{\beta\neq0}G_\beta(A(U^{n}_{\alp+\beta})-A(U^{n}_{\alp}))=c_\lambda\int_{\R^d\setminus R_0}\frac{A(\bar u(y_\alp+z,t_n))-A(\bar
u(y_\alp,t_n))}{|z|^{d+\lambda}}\ \dif z,
\end{split}
\end{equation}
where $y_\alp=x_\alp+\frac{\Delta x}{2}(1,\dots,1)$. Note that this is
an approximation of the principal value of the integral since $R_0\ra
0$ as $\Delta x\ra0$ in a symmetric way.

We now check that the numerical method \eqref{scheme} is conservative and monotone.

\begin{lemma}
The numerical method \eqref{scheme} is conservative, i.e.
\begin{align*}
\sum_{\alp\in\mathbb{Z}^d}U^{n+1}_\alp&=\sum_{\alp\in\mathbb{Z}^d}U^n_\alp\label{cons}.
\end{align*}
\end{lemma}
\begin{proof}
First we show that $\sum_{\alp\in\Z^d}
\left|U^{n}_{\alp}\right|<\infty$ for all $n=0,\dots,N$. By \eqref{scheme},
\begin{equation}\label{m7}
\begin{split}
\sum_{\alp\in\Z^d} \left|U^{n+1}_{\alp}\right|&\leq \sum_{\alp\in\Z^d} \bigg\{\left|U^{n}_{\alp}\right|+\Delta
  t\sum_{l=1}^d\left|D_l^-F_l(U_{\alp}^{n},U_{\alp+e_l}^{n})\right|\\
  &\qquad\qquad+\Delta
  t\sum_{\beta\neq0}G_\beta\left|A(U^{n}_{\alp+\beta})-A(U^{n}_{\alp})\right|\bigg\}\\
 & \leq \sum_{\alp\in\Z^d} \bigg\{\left|U^{n}_{\alp}\right|+\frac{\Delta t}{\Delta
    x}\sum_{l=1}^d\Big(L_F\left|U_{\alp}^{n}-U_{\alp-e_l}^{n}\right|+L_F\left|U_{\alp+e_l}^{n}-U_{\alp}^{n}\right|\Big)\\
& \qquad\qquad+\Delta t
\sum_{\beta\neq0}G_\beta\Big(|A(U^{n}_{\alp+\beta})|+|A(U^{n}_{\alp})|\Big)\bigg\}\\
&\leq \bigg(1+4dL_F\frac{\Delta t}{\Delta
    x}+ 2 L_A\Delta t \sum_{\beta\neq0}G_\beta
  \bigg)\sum_{\alp\in\Z^d}|U_{\alp}^{n}|,
\end{split}
\end{equation}
where, using that $\{z:|z|<\frac{\Delta x}{2}\}\subseteq R_0$,
\begin{equation*}
\begin{split}
\sum_{\beta\neq 0}G_\beta=c_\lambda\int_{\R^d\setminus R_0}\frac{\dif z}{|z|^{d+\lambda}}\leq c_\lambda\int_{|z|>\frac{\Delta
    x}{2}}\frac{\dif z}{|z|^{d+\lambda}}= c_\lambda \left(\frac{2}{\Delta
x}\right)^\lambda\int_{|z|>1}\frac{\dif z}{|z|^{d+\lambda}}.
\end{split}
\end{equation*}
Since $\Delta x^d\sum_{\alp\in\Z^d}|U_\alp^0|=\|\bar u_0\|_{L^1(\R^d)}<\infty$, we can iterate estimate \eqref{m7} to find that $\sum_{\alp\in\Z^d}
\left|U^{n}_{\alp}\right|<\infty$ and hence $\lim_{|\alp|\rightarrow\infty}|U^n_\alp|=0$ for all $n=0,\dots,N$.

Now we sum \eqref{scheme} over $\alp$ to find that
\begin{align*}
\sum_{\alp\in\mathbb{Z}^d}U^{n+1}_{\alp}=\sum_{\alp\in\mathbb{Z}^d}U^{n}_{\alp}&-\Delta
t\sum_{\alp\in\mathbb{Z}^d}\sum_{l=1}^dD_l^-F_l(U_{\alp}^{n},U_{\alp+e_l}^{n})\\
&+\Delta t\sum_{\alp\in\mathbb{Z}^d}\sum_{\beta\neq0}G_\beta(A(U^{n}_{\alp+\beta})-A(U^{n}_{\alp})).
\end{align*}
The proof is now complete if we can show that the $F$ and $G$ sums are
equal to zero. The $F$-sum is telescoping and since
$\lim_{|\alp|\rightarrow\infty}|U^n_\alp|=0$,
\begin{equation*}
\begin{split}
\sum_{\alp\in\mathbb{Z}^d}D_l^-F_l(U_{\alp}^{n},U_{\alp+e_l}^{n})=\sum_{\alp\in\mathbb{Z}^d}\frac{F(U_{\alp}^{n},U_{\alp+e_l}^{n})-F(U_{\alp-e_l}^{n},U_{\alp}^{n})}{\Delta
x}=0.
\end{split}
\end{equation*}
To treat the $G$-sum, note that we have found above that
$$\sum_{\alp\in\mathbb{Z}^d}\sum_{\beta\neq0}G_\beta
\left|A(U^{n}_{\alp+\beta})-A(U^{n}_{\alp})\right|\leq 2 L_A\Delta t \sum_{\beta\neq0}G_\beta\sum_{\alp\in\Z^d}|U_{\alp}^{n}|<\infty,$$ and we also have that
$\sum_{\alp}|A(U^{n}_{\alp})|\leq L_A\sum_{\alp}|U^{n}_{\alp}|<\infty$.  In view of this we can now change the order of summation, and split the sums to find
that
\begin{align*}
\sum_{\alp\in\mathbb{Z}^d}\sum_{\beta\neq0}G_\beta(A(U^{n}_{\alp+\beta})-A(U^{n}_{\alp}))
&=\sum_{\beta\neq0}G_\beta\sum_{\alp\in\mathbb{Z}^d}\Big(A(U^{n}_{\alp+\beta})-A(U^{n}_{\alp})\Big)\\
&=\sum_{\beta\neq0}G_\beta\Big(\sum_{\alp\in\mathbb{Z}^d}A(U^{n}_{\alp})-\sum_{\alp\in\mathbb{Z}^d}A(U^{n}_{\alp})\Big)=0.
\end{align*}
The proof is now complete.
\end{proof}
Next, we check monotonicity by showing that the right-hand side of the numerical method \eqref{scheme} is a non-decreasing function of all its variables
$U_\beta^n$.  This is clear for all $U^n_\beta$ such that $\beta\neq \alp$ since the numerical flux $F_l$ is increasing w.r.t.~ its first variable,
non-increasing w.r.t.~its second one, the function $A$ is
non-decreasing, and the weights $G_\beta$ are all positive. Then we
differentiate the right hand side of \eqref{scheme}
w.r.t.~$U^n_\alp$ and find that it is non-negative provided the
following the CFL condition holds,
\begin{equation}\label{CFL}
\begin{split}
2dL_{F}\frac{\Delta t}{\Delta x}+\Big(c_\lambda2^\lambda L_A\int_{|z|>1}\frac{\dif z}{|z|^{d+\lambda}}\Big)\frac{\Delta t}{\Delta x^{\lambda}}\leq1.
\end{split}
\end{equation}
We have thus proved the following result:
\begin{lemma}
The numerical method \eqref{scheme} is monotone provided that the CFL condition \eqref{CFL} is assumed to hold.
\end{lemma}
In what follows, the CFL condition \eqref{CFL} is always assumed to hold, and monotonicity is thus always ensured.

\subsection{Convergence toward the entropy solution}
We prove that any limit  of a uniformly bounded sequence of solutions of the
numerical method \eqref{scheme} is an entropy solution of \eqref{1}.

\begin{theorem}\label{th:conv}
If $\{\bar u\}$ is a sequence
of solutions of \eqref{scheme}, uniformly bounded in
$L^\infty(Q_T)$, and there exists $u\in
L^\infty(Q_T)\cap C([0,T];L^1(\R^d))$ such that $\bar u\ra u$ in
$C([0,T];L^1(\R^d))$ as $\Delta x\ra0$, then $u$ is an entropy
solution of  \eqref{1}.
\end{theorem}

\begin{proof}
Note that part $i)$  in the definition of entropy solution
(Definition \ref{L1-entropy}) is already satisfied. Part $iii)$
follows since $\|\bar u(\cdot,0)-u_0\|_{L^1(\R^d)}\ra0$ as
$\dx\ra0$ by the definition of $\bar u$. What remains to prove is part
$ii)$.

First we prove that the numerical method \eqref{scheme} satisfies a
discrete entropy inequality which resembles the one in $ii)$,
Definition \ref{L1-entropy}.
To this end, let us introduce the notation $a\wedge b=\min\{a,b\}$ and $a\vee b=\max\{a,b\}$, choose an $r>0$, and exploit monotonicity to obtain the
inequalities
\begin{equation*}
\begin{split}
U_{\alp}^{n+1}\vee k\leq U_{\alp}^{n}\vee k&-\Delta t\sum_{l=1}^dD_l^-F_l(U_{\alp}^{n}\vee k,U_{\alp+e_l}^{n}\vee k)\\\
&+\Delta t\sum_{0<\dx|\beta|\leq r}G_{\beta}\Big(A(U^{n}_{\alpha+\beta}\vee k)-A(U^{n}_{\alpha}\vee k)\Big)\\
&+\Delta t\, \mathbf{1}_{(k,+\infty)}(U_{\alpha}^{n+1})\sum_{\dx|\beta|>r}G_{\beta}\Big(A(U^{n}_{\alpha+\beta})-A(U^{n}_{\alpha})\Big)
\end{split}
\end{equation*}
and
\begin{equation*}
\begin{split}
U_{\alp}^{n+1}\wedge k\leq U_{\alp}^{n}\wedge k&-\Delta t\sum_{l=1}^dD_l^-F_l(U_{\alp}^{n}\wedge k,U_{\alp+e_l}^{n}\wedge k)\\\
&+\Delta t\sum_{0<\dx|\beta|\leq r}G_{\beta}\Big(A(U^{n}_{\alpha+\beta}\wedge k)-A(U^{n}_{\alpha}\wedge k)\Big)\\
&+\Delta t\, \mathbf{1}_{(-\infty,k)}(U_{\alpha}^{n+1})\sum_{\dx|\beta|>r}G_{\beta}\Big(A(U^{n}_{\alpha+\beta})-A(U^{n}_{\alpha})\Big).
\end{split}
\end{equation*}
Note that the polygonal set
$$P_r:=\underset{0<\dx|\beta|\leq r}{\bigcup}(x_\beta+R_0)$$
($x_\beta=\dx\beta$) does not include points from the box $R_0$, and
converges to the punctured ball $\{z:0<|z|\leq r\}$ as $\Delta x\ra0$
in the sense that $\mathbf{1}_{P_r}(z)\ra \mathbf{1}_{0<|z|\leq r}(z)$
a.e. as $\Delta x\ra0$.

Remember that $\eta_k(U_{\alpha}^{n})=|U_{\alpha}^{n}-k|$, and let
$$Q_{h,l}(U_{\alp}^{n})=F_l(U_{\alp}^{n}\vee k,U_{\alp+e_l}^{n}\vee k)-F_l(U_{\alp}^{n}\wedge k,U_{\alp+e_l}^{n}\wedge k).$$
Thanks to the relations
\begin{equation*}
\begin{split}
|u-k|&=u\vee k-u\wedge k,\\
|A(u)-A(k)|&=A(u\vee k)-A(u\wedge k),
\end{split}
\end{equation*}
we can subtract the above two inequalities to obtain that
\begin{equation*}
\begin{split}
\eta_k(U_{\alp}^{n+1})-\eta_k(U_{\alp}^{n})&+\frac{\Delta t}{\Delta x}\sum_{l=1}^d\Big(Q_{h,l}(U_{\alp}^{n})-Q_{h,l}(U_{\alp-e_l}^{n})\Big)\\
&-\Delta t\sum_{0<\dx|\beta|\leq r}G_{\beta}\Big(\eta_{A(k)}(A(U^{n}_{\alpha+\beta}))-\eta_{A(k)}(A(U^{n}_{\alpha}))\Big)\\
&-\Delta t\ \eta'_k(U^{n+1}_\alpha)\sum_{\dx|\beta|>r}G_{\beta}\Big(A(U^{n}_{\alpha+\beta})-A(U^{n}_{\alpha})\Big)\leq0.
\end{split}
\end{equation*}
Let us take a nonnegative function $\varphi\in C_c^\infty(Q_T)$, and define $\varphi_\alpha^n=\varphi(x_\alpha,t_n)$. If we multiply both sides of the above
inequality by $\varphi^n_\alpha$, sum over all $\alpha\in\mathbb{Z}^d$ and all $n\in\{0,\ldots,N\}$, and use summation by parts for the local terms, we end up
with the cell entropy inequality
\begin{equation}\label{disc_ent_ineq}
\begin{split}
  &\Delta x^d\Delta t\sum_{n=1}^{N}\sum_{\alp\in\mathbb{Z}^d}\eta_k(U_{\alp}^{n})\ \frac{\varphi_{\alp}^{n}-\varphi_{\alp}^{n-1}}{\Delta t}\\
  &+\Delta x^d\Delta t\sum_{n=0}^{N}\sum_{\alp\in\mathbb{Z}^d}\sum_{l=1}^dQ_{h,l}(U_\alp^n)\ \frac{\varphi_{\alp+e_l}^{n}-\varphi_{\alp}^{n}}{\Delta x}\\
  &+\Delta x^d\Delta
  t\sum_{n=0}^{N}\sum_{\alpha\in\mathbb{Z}^d}\eta_{A(k)}(A(U^{n}_{\alpha}))\sum_{0<\dx|\beta|\leq
    r}
  G_{\beta}\Big(\varphi_{\alpha+\beta}^{n}-\varphi_\alpha^{n}\Big)\\
  &+\Delta x^d\Delta t
  \sum_{n=0}^{N}\sum_{\alp\in\mathbb{Z}^d}\eta'_{k}(U^{n+1}_\alp)\
  \varphi_\alp^{n}\sum_{\dx|\beta|>r}G_\beta\Big(A(U^{n}_{\alp+\beta})-A(U^{n}_{\alp})\Big)\geq0.
\end{split}
\end{equation}
To derive this inequality we have used the
change of indices $(\beta,\alpha)\rightarrow(-\beta,\alpha+\beta)$ to
see that
\begin{equation*}
\begin{split}
&\Delta x^d\Delta t\sum_{n=0}^{N}\sum_{\alpha\in\mathbb{Z}^d}\varphi_\alpha^{n}\sum_{0<\dx|\beta|\leq r}
G_{\beta}\Big(\eta_{A(k)}(A(U^{n}_{\alpha+\beta}))-\eta_{A(k)}(A(U^{n}_{\alpha}))\Big)\\
&=\Delta x^d\Delta t\sum_{n=0}^{N}\sum_{\alpha\in\mathbb{Z}^d}\eta_{A(k)}(A(U^{n}_{\alpha}))\sum_{0<\dx|\beta|\leq r}
G_{\beta}\Big(\varphi_{\alpha+\beta}^{n}-\varphi_\alpha^{n}\Big).
\end{split}
\end{equation*}

Let $R_\alpha=x_\alpha+R$. We now claim that for each fixed $\Delta x>0$, inequality \eqref{disc_ent_ineq} implies
\begin{equation}\label{eeemmm}
\begin{split}
&\sum_{n=0}^{N-1}\sum_{\alp\in\mathbb{Z}^d}\eta_k(U_{\alp}^{n})\int_{t_n}^{t_{n+1}}\int_{R_\alpha}\frac{\varphi(x,t)-\varphi(x,t-\Delta t)}{\Delta t}\ \dif x\dif t\\
&+\sum_{n=0}^{N-1}\sum_{\alp\in\mathbb{Z}^d}\sum_{l=1}^dQ_{h,l}(U_\alp^n)\int_{t_n}^{t_{n+1}}\int_{R_\alpha}
\frac{\varphi(x+\Delta x\, e_l,t)-\varphi(x,t)}{\Delta x}\ \dif x\dif t\\
&+\sum_{n=0}^{N-1}\sum_{\alpha\in\mathbb{Z}^d}\eta_{A(k)}(A(U^{n}_{\alpha}))\int_{t_n}^{t_{n+1}}\int_{R_\alpha}
\int_{P_r}\frac{\varphi(x+z,t)-\varphi(x,t)}{|z|^{d+\lambda}}\ \dif z\dif x\dif t\\
&+\sum_{n=0}^{N-1}\sum_{\alp\in\mathbb{Z}^d}\eta'_{k}(U^{n+1}_\alp)\sum_{\dx|\beta|>r}G_\beta\Big(A(U^{n}_{\alp+\beta})-A(U^{n}_{\alp})\Big)
\int_{t_n}^{t_{n+1}}\int_{R_\alpha}\varphi(x,t)\ \dif x\dif t\geq0.
\end{split}
\end{equation}
To see this we proceed by contradiction, and assume that \eqref{eeemmm} is strictly negative. We then sum together several inequalities of the form
\eqref{disc_ent_ineq} where, instead of $\varphi_\alpha^n=\varphi(x_\alpha,t_n)$ which are computed on the original space/time grid $(x_\alpha,t_n)$, we use
the values $\varphi^n_\alpha=\varphi(\hat x_\alpha,\hat t_n)$ computed on the finer grid $(\hat x_\alpha,\hat t_n)$ where $\hat x_\alpha=(\Delta x/M)\alpha$
while $\hat t_n=n(\Delta t/M)$ for some $M>0$. Note that, since all these inequalities of the form \eqref{disc_ent_ineq} share the same underlying numerical
solution $(U_i^n)$, they can be rearranged as one inequality, i.e.
\begin{equation}\label{j3}
\begin{split}
&\sum_{n=1}^{d}\sum_{\alp\in\mathbb{Z}^d}\eta_k(U_{\alp}^{n})\left(\left(\frac{\Delta x}{M}\right)^d\frac{\Delta t}{M}
\sum_{m:\, \hat t_m\in[t_n,t_{n+1})}\sum_{\gamma:\, \hat x_\gamma\in R_\alpha}\frac{\varphi_{\gamma}^{m}-\varphi_{\gamma}^{m-1}}{\Delta t}\right)\\
&+\sum_{n=0}^{d}\sum_{\alp\in\mathbb{Z}^d}\sum_{l=1}^dQ_{h,l}(U_\alp^n)\left(\left(\frac{\Delta x}{M}\right)^d\frac{\Delta t}{M}
\sum_{m:\, \hat t_m\in[t_n,t_{n+1})}\sum_{\gamma:\, \hat x_\gamma\in R_\alpha}\frac{\varphi_{\gamma+e_l}^{m}-\varphi_{\gamma}^{m}}{\Delta x}\right)\\
&+\sum_{n=0}^{d}\sum_{\alpha\in\mathbb{Z}^d}\eta_{A(k)}(A(U^{n}_{\alpha}))\\
&\qquad\qquad\qquad\left(\left(\frac{\Delta x}{M}\right)^d\frac{\Delta t}{M}\sum_{m:\, \hat t_m\in[t_n,t_{n+1})}\sum_{\gamma:\, \hat x_\gamma\in
R_\alpha}\sum_{0<\dx|\beta|\leq r}
G_{\beta}\Big(\varphi_{\gamma+\beta}^{m}-\varphi_\gamma^{m}\Big)\right)\\
&+\sum_{n=0}^{d}\sum_{\alp\in\mathbb{Z}^d}\eta'_{k}(U^{n+1}_\alp)\sum_{\dx|\beta|>r}G_\beta\Big(A(U^{n}_{\alp+\beta})-A(U^{n}_{\alp})\Big)\\
&\qquad\qquad\qquad\qquad\qquad\qquad\qquad\qquad\left(\left(\frac{\Delta x}{M}\right)^d\frac{\Delta t}{M} \sum_{m:\, \hat t_m\in[t_n,t_{n+1})}\sum_{\gamma:\,
\hat x_\gamma\in R_\alpha}\varphi_\gamma^{m}\right)\geq0
\end{split}
\end{equation}
(loosely speaking, by summing all these inequalities of the form \eqref{disc_ent_ineq} together we are filling the mesh-sets $R_\alpha\times [t_n,t_{n+1})$
with several samples of the test function $\varphi$; this has been done in order to recreate in each mesh-set a Riemann sum approximation which gets closer and
closer to its respective integral as the value of the control parameter $M$ increases). The Riemann sum approximations in the first, second, and fourth term of
\eqref{j3} are arbitrarily close to their respective terms in \eqref{eeemmm} as $M$ increases. For the third term in \eqref{j3} note that,
cf.~\eqref{scheme_fund_rel},
\begin{equation}\label{a5}
\begin{split}
&\left(\frac{\Delta x}{M}\right)^d\frac{\Delta t}{M}\sum_{m:\, \hat t_m\in[t_n,t_{n+1})}\sum_{\gamma:\, \hat x_\gamma\in
R_\alpha}\sum_{0<\dx|\beta|\leq r} G_{\beta}\Big(\varphi_{\gamma+\beta}^{m}-\varphi_\gamma^{m}\Big)\\
&=\left(\frac{\Delta x}{M}\right)^d\frac{\Delta t}{M}\sum_{m:\, \hat t_m\in[t_n,t_{n+1})}\sum_{\gamma:\, \hat x_\gamma\in R_\alpha}\int_{z\in
P_r}\frac{\bar\varphi(y_\gamma+z,\hat t_m)-\bar\varphi(y_\gamma,\hat t_m)}{|z|^{d+\lambda}}\ \dif z
\end{split}
\end{equation}
(the definitions of $\bar\varphi,y_\gamma$ are analogous to those of $\bar u,y_\alpha$) and so the Riemann sum approximation on the right-hand side of
\eqref{a5} is, as $M$ increases, arbitrarily close to
\begin{equation*}
\begin{split}
\int_{t_n}^{t_{n+1}}\int_{R_\alpha}\int_{z\in P_r}\frac{\varphi(x+z,t)-\varphi(x,t)}{|z|^{d+\lambda}}\ \dif z\dif x\dif t.
\end{split}
\end{equation*}
This is due to the fact that, since we are integrating away from the singularity, the right-hand side of \eqref{a5} is well defined, and the sum over all
$(\hat x_\gamma,\hat t_m)$ can be moved inside the integral $z\in P_r$. Therefore, since \eqref{j3} is arbitrarily close to the left-hand side of
\eqref{eeemmm}, the left-hand side of \eqref{eeemmm} cannot be negative, and we have produced a contradiction.

Using the piecewise constant space/time interpolation $\bar u$, we can now rewrite inequality \eqref{eeemmm}  as
\begin{equation}\label{a4}
\begin{split}
&\iint_{Q_T}\Bigg\{\eta_k(\bar u(x,t))\, \partial_t\varphi(x,t)+\sum_{l=1}^dQ_{h,l}(\bar u(x,t))\, \partial_{x_l}\varphi(x,t)\\
&+\eta_{A(k)}(A(\bar u(x,t)))\int_{P_r}\frac{\varphi(x+z,t)-\varphi(x,t)}{|z|^{d+\lambda}}\ \dif z\\
&+\eta_{k}'(\bar u(x,t+\Delta t))\, \varphi(x,t)\int_{\R^d\setminus P_r}\frac{A(\bar u(x+z,t))-A(\bar
  u(x,t))}{|z|^{d+\lambda}}\ \dif z\Bigg\}\ \dif x\dif t\\
&\geq O(\Delta x)+O(\Delta t).
\end{split}
\end{equation}
Convergence up to a subsequence for the first three terms in
\eqref{a4} is immediate thanks to the a.e.~convergence of $\bar u$
toward $u$. For the local terms this is already well known,
cf.~\cite[Theorem 3.9]{Holden/Risebro}. For the term containing
the inner integral $P_r$, convergence follows thanks to the
convergence of $\mathbf 1_{
  P_r}\ra \mathbf 1_{0<|z|<r}$ a.e., the
properties of $\varphi$ ($\int_{P_r}\frac{\varphi(x+z,t)-\varphi(x,t)}{|z|^{d+\lambda}}\ \dif
 z$ is uniformly bounded and compactly supported),
uniform boundedness of $\bar
u$, and the fact that the function $\eta_{k}(\cdot)$ is continuous.

To conclude, we need to establish convergence for the term containing the discontinuous sign function $\eta'_k(\cdot)$, and we argue as in \cite{Droniou}
(p.~109). First note that since $\bar u\ra u$ in $C([0,T];L^1(\R^d))$,
also $\bar u(\cdot,\cdot+\Delta t)\ra u$ in $C([0,T];L^1(\R^d))$ and
a.e. for a subsequence. Then note that $\eta_{k}'(s)$ is continuous
for $s\neq k$, and that the measure of the set
$$\mathcal{U}_k=\{(x,t)\in Q_T:\ u(x,t)=k\}$$
is $0$ for a.e.~$k\in\mathbb{R}$. For such $k$, $\eta_{k}'(\bar
u(\cdot+\Delta t,\cdot))\ra \eta_{k}'(u)$ a.e., and we can go to the
limit in the term involving $\eta_k'$ in \eqref{a4} using the
dominated convergence theorem, $|\eta_k'|\leq 1$, and uniform
boundedness of $\bar u$ and $A(\bar u)$.

 For the
remaining $k$, we use an approximating sequence made of those $k$ for
which convergence holds true. To be more precise, let $a_m,b_m$ be
sequence of values such that
$\text{meas}(\mathcal{U}_{a_m})=\text{meas}(\mathcal{U}_{b_m})=0$,
where $a_m\nearrow k$ and $b_m\searrow k$. Note that the mean value
$$\frac{1}{2}(\eta'_{a_m}(u)+\eta'_{b_m}(u))\rightarrow \eta'_k(u)\qquad\text{ as $a_m,b_m\rightarrow k$}.$$
Thus we can use the entropy inequality for the sequence $a_m$ and the
entropy inequality for the sequence $b_m$, take the average, and go to
the limit to prove the entropy inequality for every critical value
$k$. Convergence for the whole sequence $\bar u$ is a consequence of
uniqueness for entropy solutions of \eqref{1}.
\end{proof}

\subsection{BV initial data: Compactness and existence.}
We now show that the sequence of solutions of the method, $\{\bar u:\Delta x>0\}$, is relatively compact whenever
$$u_0\in L^\infty(\R^d)\cap L^1(\R^d)\cap BV(\R^d).$$
Using this result and Theorem \ref{th:conv}, we then obtain existence of an entropy solution of the initial value problem \eqref{1}. We start by the following
a priori estimates.
\begin{lemma}\label{lem:apriorest}
If $u_{0}\in L^{\infty}(\mathbb{R}^d)\cap L^{1}(\mathbb{R}^d)\cap BV(\mathbb{R}^d)$, then, for all $t,s\geq 0$,
\begin{itemize}
\item[\emph{i)}]$\|\bar{u}(\cdot,t)\|_{L^{\infty}(\mathbb{R}^d)}\leq\|u_{0}\|_{L^{\infty}(\mathbb{R}^d)}$,
\item[\emph{ii)}]$\|\bar{u}(\cdot,t)\|_{L^{1}(\mathbb{R}^d)}\leq\|u_{0}\|_{L^{1}(\mathbb{R}^d)}$,
\item[\emph{iii)}]$|\bar{u}(\cdot,t)|_{BV(\mathbb{R}^d)}\leq|u_{0}|_{BV(\mathbb{R}^d)}$,
\item[\emph{iv)}]$\|\bar{u}(\cdot,s)-\bar{u}(\cdot,t)\|_{L^{1}(\mathbb{R}^d)}\leq\sigma(|s-t|+\Delta
  t)$ where, for some $c>0$,
$$\sigma(s)=\begin{cases}c\, |s|&\text{if }\lambda\in(0,1),\\
c\, |s\ln s|&\text{if }\lambda=1,\\
c\, |s|^{\frac1{\lambda}}&\text{if }\lambda\in(1,2).\end{cases}$$
\end{itemize}
\end{lemma}

Lemma \ref{lem:apriorest} along with a Kolmogorov type of compactness
theorem, cf. Theorem A.8 in ~\cite{Holden/Risebro},
yields the existence of a subsequence $\{\bar{u}\}$ which converges  in $C([0,T];L^1_{\mathrm{loc}}(\R^d))$ (and hence a.e.~up
to a further subsequence) toward
a limit $u$ as
$\Delta x\ra0$. Moreover, the limit $u$ inherits all the a priori
estimates \emph{i)-iv)} in Lemma \ref{lem:apriorest} (with $\Delta
t=0$). Moreover, by $ii)$ and the dominated convergence theorem, we see that
$\bar u\ra u$ also in $C([0,T];L^1(\R^d))$. In short, we have the
following result:
\begin{lemma}\label{lem:comp}
The numerical solutions $\{\bar u:\Delta x>0\}$ converge, up to a subsequence, toward a limit $u$ in $C([0,T];L^1(\R^d))$ as $\Delta x\rightarrow 0$. Moreover,
$$u\in L^\infty(Q_T)\cap C([0,T];L^{1}(\mathbb{R}^d))\cap L^\infty(0,T;BV(\mathbb{R}^d)).$$
\end{lemma}

Lemma \ref{lem:comp} and Theorem \ref{th:conv} imply the following
existence result:

\begin{theorem}\emph{(Existence)}\label{th:existence}
If $u_0\in L^{\infty}(\mathbb{R}^d)\cap L^{1}(\mathbb{R}^d)$, then
there exists an entropy solution of the initial value problem
\eqref{1}.
\end{theorem}

\begin{proof}
For initial data $u_0\in L^{\infty}(\mathbb{R}^d)\cap
L^{1}(\mathbb{R}^d)\cap BV(\mathbb{R}^d)$ existence is granted by the
numerical method \eqref{scheme} (Lemma \ref{lem:comp}). For more general initial data $u_0\in
L^{\infty}(\mathbb{R}^d)\cap L^{1}(\mathbb{R}^d)$, we consider approximations $u_{0,n}\in L^{\infty}(\mathbb{R}^d)\cap L^{1}(\mathbb{R}^d)\cap BV(\mathbb{R}^d)$ such that
\begin{equation*}
\begin{split}
\|u_0-u_{0,n}\|_{L^1(\R^d)}\rightarrow 0\text{ as $n\rightarrow\infty$}.
\end{split}
\end{equation*}
Let $u_m,u_n$ denote the entropy solutions corresponding to $u_{0,n},u_{0,m}$
respectively, and use the $L^1$-contraction (Theorem
\ref{L1contraction}) to see that
\begin{equation*}
\begin{split}
\|u_n-u_m\|_{C([0,T];L^1(\R^d))}\leq\|u_{0,n}-u_{0,m}\|_{L^1(\R^d)}\rightarrow 0\text{ as $n,m\rightarrow\infty$}.
\end{split}
\end{equation*}
Therefore, the sequence of entropy solutions $\{u_n\}$ is Cauchy in $C([0,T];L^1(\R^d))$ and admits a limit $u$. To prove that $u$ is also an entropy solution of
\eqref{1}, one can pass to the limit $n\rightarrow \infty$ in the entropy inequality for $u_n$.
\end{proof}

\begin{proof}[Proof of Lemma \ref{lem:apriorest}]
The maximum principle \emph{i)} is a direct consequence of monotonicity. To see this let $s=\sup_{\alpha\in\Z^d}|U_\alpha^n|$, and choose $U^n\equiv s$ to
obtain that
\begin{equation*}
\begin{split}
U_{\alp}^{n+1}\leq s-\Delta t\sum_{l=1}^dD^-_lF_l(s,s)+\Delta t\sum_{\beta\neq 0}G_\beta(A(s)-A(s))=s.
\end{split}
\end{equation*}
Similarly, choosing $U^n\equiv -s$, one obtains $U_{\alp}^{n+1}\geq -s$. Furthermore, since the numerical method \eqref{scheme} is conservative, monotone, and
translation invariant (translation invariance is a consequence of the fact that the numerical method does not explicitly depend on the variables $x_\alp,t_n$),
inequalities \emph{ii)-iii)} are consequences of the results due to Crandall-Tartar \cite{Crandall/Tartar}.

We now prove \emph{iv)}. By \eqref{scheme} and Lipschitz continuity of $F_l$,
\begin{equation*}
\begin{split}
&\left|U_{\alp}^{n+1}-U_{\alp}^{n}\right|\\
&\leq\frac{\Delta tL_F}{\Delta x}\sum_{l=1}^d\bigg(|U^n_{\alp+e_l}-U^n_{\alp}|+|U^n_{\alp}-U^n_{\alp-e_l}|\bigg)+\Delta
t\sum_{\beta\neq0}G_\beta|A(U^{n}_{\alp+\beta})-A(U^{n}_{\alp})|.
\end{split}
\end{equation*}
Let us multiply by $\Delta x^d$ in the above inequality, and sum over $\alp\in\Z^d$ to see that
\begin{align*}
&\Delta x^d\sum_{\alp\in\Z^d}\left|U_{\alp}^{n+1}-U_{\alp}^{n}\right|\\
&\leq 2L_F\Delta x^{d-1}\Delta t\sum_{l=1}^d\sum_{\alp\in\Z^d}|U^n_{\alp+e_l}-U^n_{\alp}|+\Delta x^d\Delta
t\sum_{\alp\in\Z^d}\sum_{\beta\neq0}G_\beta|A(U^{n}_{\alp+\beta})-A(U^{n}_{\alp})|.
\end{align*}
Let $\bar u^n(\cdot)=\bar u(\cdot,t_n)$ and note that the first term then
is equal to
$$2L_F\Delta t \sum_{l=1}^d\int_{\R^{d-1}}|\bar
u^n(\cdot,x')|_{BV_{x_l}(\R)}\ \dif x'\leq  2dL_F\Delta t|\bar u^n|_{BV(\R^d)}=O(\Delta t),$$ while the second term can be estimated by
(cf.~\eqref{scheme_fund_rel})
\begin{equation*}
\begin{split}
&\Delta x^d\sum_{\alp\in\mathbb{Z}^d}\sum_{\beta\neq0}G_\beta\left|A(U^{n}_{\alp+\beta})-A(U^{n}_{\alp})\right|\\
&\leq c_\lambda\sum_{\alp\in\mathbb{Z}^d}\int_{|z|>\frac{\Delta x}{2}}\frac{|A(\bar
  u^n(y_\alp+z))-A(\bar u^n(y_\alp))|}{|z|^{d+\lambda}}\ \dif z\Delta x^d\\
&\leq c_\lambda L_A\bigg(\int_{\frac{\Delta
    x}{2}<|z|<1}+\int_{|z|>1}\bigg)\sum_{\alp\in\mathbb{Z}^d}\frac{|\bar
  u^n(y_\alp+z)-\bar u^n(y_\alp)|}{|z|^{d+\lambda}}\ \Delta x^d\dif z\\
&\leq c_\lambda L_A\bigg(|\bar u^n|_{BV(\mathbb{R}^d)}\int_{\frac{\Delta x}{2}<|z|<1}\frac{|z|}{|z|^{d+\lambda}}\ \dif z+2\|\bar
u^n\|_{L^1(\mathbb{R}^d)}\int_{|z|>1}\frac{\dif z}{|z|^{d+\lambda}}\bigg).
\end{split}
\end{equation*}
Easy computations in polar coordinates show that the second integral is $O(1)$ while
$$I_{\Delta x}=\int_{\frac{\Delta
    x}{2}<|z|<1}\frac{|z|}{|z|^{d+\lambda}}\ \dif z=
\begin{cases}
O(1)&\text{if }\lambda\in(0,1),\\
O(|\ln\Delta x|)&\text{if }\lambda=1,\\
O(\Delta x^{1-\lambda})&\text{if }\lambda\in(1,2).
\end{cases}$$
Adding all the above estimates yields
\begin{align*}
\Delta x^d\sum_{\alp\in\Z^d}|U_{\alp}^{n+1}-U_{\alp}^{n}|=O(\Delta t)+ O(\Delta t I_{\Delta x})+O(\Delta t).
\end{align*}
By the CFL condition \eqref{CFL}, $\Delta t I_{\Delta
  x}=\sigma(\Delta t)$, and the result follows.
\end{proof}

%%%%%%%%%%%%%%%%%%%%%%%%%%%%%%%%%%%%%%%%%%%%%%%%%%%%%%%%%%%%%%%%%%%%%%%%%%%%%%%%%%%%%%%%%%%%%%%%%%%%%%%%%%%%%%%%%%%%%%%%%%%%%%%%%%%%%%%%%%%%%%%%%%%%%%%%%%%%%%%%%%%%%%%%%%%
%%%%%%%%%%%%%%%%%%%%%%%%%%%%%%%%%%%%%%%%%%%%%%%%%%%%%%%%%%%%%%%%%%%%%%%%%%%%%%%%%%%%%%%%%%%%%%%%%%%%%%%%%%%%%%%%%%%%%%%%%%%%%%%%%%%%%%%%%%%%%%%%%%%%%%%%%%%%%%%%%%%%%%%%%%%

\section{Extension to general L\'{e}vy operators}
\label{sec:Levy}
The ideas developed in this paper can also be used to establish well-posedness for entropy solutions of a more general class of fractional equations of the
form
\begin{align}\label{1bis}
\begin{cases}
\partial_tu+\nabla\cdot f(u)=g_\mu[A(u)]&\text{in }Q_T=\mathbb{R}^d\times(0,T),\\
u(x,0)=u_{0}(x)&\text{in }\mathbb{R}^d,
\end{cases}
\end{align}
where the fractional Laplacian $g$ has been replaced with a more general L\'{e}vy operator $g_\mu$:
$$g_\mu[\phi](x)=\int_{|z|>0}\phi(x+z)-\phi(x)-z\cdot
\nabla\phi(x)\mathbf{1}_{|z|<1}\ \dif\mu(z),$$ where the L\'{e}vy measure $\mu$ is a positive Radon measure satisfying
\begin{align}\label{ooo}
\int_{|z|>0}|z|^2\wedge 1 \ \mu(\dif z) <\infty.
\end{align}
Note that $g_\mu$ is self-adjoint if and only if $\mu$ is symmetric:
$\mu(-B)=\mu(B)$ for all open sets $B$. The
adjoint $g_\mu^*$ (defined through $\int u g_\mu[v] = \int g_\mu^*[u] v$) equals
$$g_\mu^*[\phi](x)=\int_{|z|>0}\phi(x-z)-\phi(x)+z\cdot
\nabla\phi(x)\mathbf{1}_{|z|<1}\ \dif\mu(z).$$ A Taylor expansion shows that both $g_\mu[\phi]$ and $g_\mu^*[\phi]$ are well defined whenever $\phi$ is $C^2$
and bounded. The gradient term is needed when $\mu$ is not radially symmetric, in the radially symmetric case $g_\mu$ $(=g_\mu^*)$ can be defined as before as
a principal value and no gradient term. The operator $g_\mu$ is the generator of a pure jump L\'{e}vy process and these processes have many applications in
Physics and Finance, cf.~e.g.~\cite{Ap:Book}.

We need a modified definition of Entropy solutions.
Remember the notation $\eta_k$ and $q_k$ introduced in
Section \ref{sec:def}, and define for $r>0$,
$$g_\mu[\varphi]=g_{\mu,r}[\varphi]+g^r_\mu[\varphi]-\gamma^r_\mu\cdot
\nabla\varphi$$ where $g_{\mu,r}[\varphi](x)=g_\mu[\varphi(\cdot)
\mathbf{1}_{|z|\leq r}](x)$,
\begin{align*}
g^r_\mu[\varphi](x)=\int_{|z|>r}\varphi(x+z)-\varphi(x)\ \mu(\dif z),\quad\text{and}\quad \gamma^r_\mu=\int_{r<|z|<1} z\ \mu(\dif z).
\end{align*}
We also use the notation $g_{\mu,r}^*$ and $g^{r,*}_\mu$ for the
adjoint operators, and note that
$$g_\mu^*[\phi]=g_{\mu,r}^*[\phi]+g^{r,*}_\mu[\phi]+\gamma^r_\mu\cdot
\nabla\phi.$$ Let us point out that the adjoint operator $g_{\mu}^*$ could have also been defined as $g_{\nu}$ with $\nu(B)=\mu(-B)$. From this equivalent
definition it is clear that the adjoint operator $g_{\mu}^*$ is still a L\'{e}vy operator.

\begin{definition}
A function $u$ is an entropy solution of the initial value problem \eqref{1bis} provided that
\begin{itemize}
\item[\emph{i)}]$u\in L^{\infty}(Q_T)\cap C([0,T];L^{1}(\mathbb{R}^d))$;
\item[\emph{ii)}]for all $k\in\mathbb{R}$, all $r>0$, and all nonnegative test functions $\varphi\in C_{c}^{\infty}(Q_T)$,
\begin{equation*}
\begin{split}\label{lll}
&\iint_{Q_T}\eta_k(u)\partial_t\varphi+q_k(u)\cdot\nabla\varphi+\eta_{A(k)}(A(u))\,
g_{\mu,r}^*[\varphi]\\
&\qquad+\eta'_k(u)\, g^r_\mu[A(u)]\, \varphi+\eta_{A(k)}(A(u))\; \gamma_\mu^r\cdot \nabla\varphi\ \dif x\dif t\geq0;
\end{split}
\end{equation*}
\item[\emph{iii)}] $u(x,0)=u_0(x)$ a.e.
%$\mathrm{esslim}_{t\ra0}\|u(\cdot,t)-u_0(\cdot)\|_{L^1(\mathbb{R}^d)}=0$.
\end{itemize}
\end{definition}
\begin{remark}
All terms in $ii)$ are well-defined in view of $i)$. Except for the $g^r_\mu$-term, this follows from the discussion proceeding Definition \ref{L1-entropy} --
see Remark \ref{WDef}. Note that the integrand of $g^r_\mu[A(u)]$ is measurable w.r.t.~the product measure $\dif \mu(z)\dif x\dif t$ since since it is the
$\dif \mu(z)\dif x\dif t$-a.e.~limit of continuous functions. This follows readily from the fact that $u$ is the $\dif x\dif t$-a.e.~limit of smooth functions.
Integrability then follows by Fubini's theorem, integrate first w.r.t.~to $\dif x\dif t$ and then w.r.t.~$\dif \mu(z)$ using \eqref{ooo}. By Fubini we also see
that $g^r_\mu[A(u)]\in C([0,T];L^1(\R^d))$ and it easily follows that the $g^r_\mu$-term is well-defined.
\end{remark}

Again classical solutions are entropy solutions and entropy solutions are weak solutions. The proof is essentially the same as the one given in Section
\ref{sec:def} with the additional information that whenever $A(u)$ is smooth
$$\eta_k'(u(x))\nabla[A(u(x))]=\eta_{A(k)}'(A(u(x))\nabla[A(u(x))]=\nabla
[\eta_{A(k)}(A(u(x)))]\ a.e.$$
We also have a $L^1$-contraction and hence uniqueness result:
\begin{theorem}
  Let $u$ and $v$ be two entropy solutions of the initial value
  problem \eqref{1bis} with initial data $u_0$ and $v_0$. Then, for
  all $t\in(0,T)$,
\begin{equation*}
\|u(\cdot,t)-v(\cdot,t)\|_{L^1(\mathbb{R}^d)}\leq\|u_0-v_0\|_{L^1(\mathbb{R}^d)}.
\end{equation*}
\end{theorem}

\begin{proof}
  We proceed as in the proof of Theorem \ref{L1contraction}: let us
  take the entropy inequality for $u=u(x,t)$ and the one for
  $v=v(y,s)$, integrate both in space/time, and sum the resulting
  inequalities together to obtain an expression equivalent to
  \eqref{julie6}. At this point we use the change of variables
  $(x,y)\rightarrow(x-z,y-z)$ to obtain the inequality
\begin{equation*}
\begin{split}
\iint_{Q_T}\iint_{Q_T}&\eta(u(x,t),v(y,s))\, (\partial_t+\partial_s)\psi(x,y,t,s)\\
&+q(u(x,t),v(y,s))\cdot(\nabla_x+\nabla_y)\psi(x,y,t,s)\\
&+\eta(A(u(x,t)),A(v(y,s)))\ g_{\mu,r}^\ast[\psi(\cdot,y,t,s)](x)\\
&+\eta(A(u(x,t)),A(v(y,s)))\ g_{\mu,r}^\ast[\psi(x,\cdot,t,s)](y)\\
&+\eta(A(u(x,t)),A(v(y,s)))\ \gamma_\mu^r\cdot (\nabla_x+\nabla_y)\psi(x,y,t,s)\\
&+\eta(A(u(x,t)),A(v(y,s)))\ \tilde g_\mu^{r,\ast}[\psi(\cdot,\cdot,t,s)](x,y)\ \dif w\geq0,
\end{split}
\end{equation*}
where
\begin{equation*}
\begin{split}
\tg_\mu^{r,\ast}[\varphi(\cdot,\cdot)](x,y)=\int_{|z|>r}\varphi(x-z,y-z)-\varphi(x,y)\ \mu(\dif z).
\end{split}
\end{equation*}
We can now send $r\rightarrow 0$ and recover the equivalent of expression \eqref{ineq} in the present setting,
\begin{equation*}
\begin{split}
\iint_{Q_T}\iint_{Q_T}&\eta(u(x,t),v(y,s))\, (\partial_t+\partial_s)\psi(x,y,t,s)\\
&+q(u(x,t),v(y,s))\cdot(\nabla_x+\nabla_y)\psi(x,y,t,s)\\
&+\eta(A(u(x,t)),A(v(y,s)))\ \tilde g_\mu^\ast[\psi(\cdot,\cdot,t,s)](x,y)\ \dif w\geq0,
\end{split}
\end{equation*}
where
\begin{equation*}
\begin{split}
&\tg_\mu^\ast[\varphi(\cdot,\cdot)](x,y)\\
&=\int_{|z|>0}\varphi(x-z,y-z)-\varphi(x,y)+z\cdot (\nabla_x+\nabla_y)\varphi(x,y)\mathbf{1}_{|z|<1}\ \mu(\dif z).
\end{split}
\end{equation*}
 From now on, the proof follows the one of Theorem \ref{L1contraction}
 (just replace the operator $g$ therein with the operator
 $g_\mu^\ast$).
\end{proof}

Existence of solutions can be obtained e.g.~by the vanishing viscosity method and a compensated compactness argument, but we do not give the details here. We
just remark that the vanishing viscosity equations have smooth solutions since the principle term is the (linear 2nd order) Laplace term.

\begin{theorem}
There exists a unique entropy solution of the initial value problem
\eqref{1bis}.
\end{theorem}

\section{Connections to HJB equations}
\label{sec:HJB} In one space dimension it is well known that the gradient of the (viscosity) solution of a HJB equation is the (entropy) solution of a
conservation law, see e.g.~\cite{Li:Book}. Variants of this result are still true in the current fractional setting as we will explain now. First we consider
the following two initial value problems in one space dimension:
\begin{equation}
\tag{HJB}
\left\{
\begin{split}
u_t + f(\del_xu)+g[u]&=\eps\del_x^2u \quad\  \text{in } Q_T,\\
u(x,0)&=u_0(x) \quad \text{in } \R,
\end{split}
\right.
\end{equation}
and
\begin{equation}
\tag{FCL}
\left\{
\begin{split}
v_t + \del_xf(v)+g[v]&=\eps\del_x^2v \qquad\ \text{in } Q_T,\\
v(x,0)&=\del_xu_0(x) \quad \text{in } \R,
\end{split}
\right.
\end{equation}
for any $\eps\geq0$. The first equation
is a HJB equation and the second one a fractional conservation
law. To simplify, let us consider the following strong but rather
standard regularity assumptions:
\begin{itemize}
\item[(a1)] $f\in C^2(\R)$,
\item[(a2)] $u_0$ is bounded and Lipschitz continuous, and
\item[(a3)] $\del_xu_0$ is bounded and belongs to $L^1(\R)\cap BV(\R)$.
\end{itemize}
Standard results then show that:
\begin{itemize}
\item[(i)] there is a unique bounded H\"older continuous (viscosity) solution
  $u^\eps$ of (HJB) for any $\eps\geq0$ \cite{JK:ContDep},
\item[(ii)] there is a unique bounded (entropy) solution
  $v^\eps\in L^1(0,T;BV\cap L^1)$ of the (fractional) conservation law
  for any $\eps\geq0$ \cite{Alibaud,Droniou/Imbert},
\item[(iii)] when $\eps>0$ both $u_\eps$ and $v_\eps$ are $C^2$,
\item[(iv)] $u^\eps\ra u^0$ uniformly \cite{JK:ContDep} and $v^\eps\ra
  v^0$ in $L^1$ \cite{Alibaud} as $\eps\ra0$.
\end{itemize}
By differentiating (HJB) and using uniqueness for (FCL), we find that
$$v^\eps=\del_xu^\eps$$
for any $\eps>0$, and hence
$$\iint v^\eps\phi=-\iint u^\eps\del_x\phi \quad\text{for
  any}\quad\phi\in C_c^\infty(Q_T).$$
Sending $\eps\ra 0$ in the above inequality using dominated
convergence theorem and (iv) then leads to
$$\iint v^0\phi=-\iint u^0\del_x\phi \quad\text{for
  any}\quad\phi\in C_c^\infty(Q_T),$$
and we have the following result:
\begin{theorem}
The distributional $x$-derivative of the viscosity solution of (HJB)
is equal to the unique entropy solution of (FCL).
\end{theorem}

The only part missing in the proof of this theorem, is the proof of (iii). This result follows e.g.~from energy estimates and standard parabolic compactness
results (yields $L^2(0,T;H^1)$ solutions) combined with regularity theory  for the Heat equation, interpolation, and bootstrapping arguments (yields smooth
solutions). We skip the long and fairly standard details.

If we drop the convection term, we get a similar correspondence in
any space dimension. Consider the following two initial value problems:
\begin{equation*}
\tag{HJB2}
\left\{
\begin{split}
u_t - A(g_\mu[u])&=\eps\Delta u \quad\ \, \text{in } Q_T,\\
u(x,0)&=u_0(x) \quad \text{in } \R^d,
\end{split}
\right.
\end{equation*}
and
\begin{equation*}
\tag{FDE}
\left\{
\begin{split}
v_t -g_\mu[A(u)]&=\eps\Delta v \qquad\quad \text{in } Q_T,\\
v(x,0)&=g_\mu[u_0](x) \quad \text{in } \R^d,
\end{split}
\right.
\end{equation*}
for any $\eps\geq0$. The first equation is still a HJB equation while the second one is a degenerate fractional diffusion equation. To simplify, let us
consider the following rather strong regularity assumptions:
\begin{itemize}
\item[(b1)] $A\in C^2(\R)$ is non-decreasing and Lipschitz continuous,
\item[(b2)] $u_0$ is bounded and Lipschitz continuous, and
\item[(b3)] $g_\mu[u_0]$ is bounded, BV, and belongs to $L^1$.
\end{itemize}
Again we have
the following of properties:
\begin{itemize}
\item[(i)] there is a unique bounded H\"older continuous (viscosity) solution
  $u^\eps$ of (HJB2) for any $\eps\geq0$,
\item[(ii)] there is a unique bounded (entropy) solution
  $v^\eps\in L^1(0,T;BV\cap L^1)$ of (FDE) for any $\eps\geq0$,
\item[(iii)] when $\eps>0$ both $u_\eps$ and $v_\eps$ are $C^2$,
\item[(iv)] $u^\eps\ra u^0$ uniformly and $v^\eps\ra v^0$ in $L^1$ as
  $\eps\ra0$.
\end{itemize}
By applying $g_\mu$ to (HJB2) and using uniqueness for (FDE), we find that
$$v^\eps=g_\mu[u^\eps]$$
for any $\eps>0$, and hence
$$\iint v^\eps\phi=\iint u^\eps g_\mu^*[\phi] \quad\text{for
  any}\quad\phi\in C_c^\infty(Q_T).$$
Sending $\eps\ra 0$ in the above inequality using the dominated
convergence theorem and (iv) then leads to
$$\iint v^0\phi=\iint u^0g_\mu^*[\phi]  \quad\text{for
  any}\quad\phi\in C_c^\infty(Q_T),$$
and we have the following result:
\begin{theorem}
If $u$ is the unique viscosity solution of (HJB2), then
$v=g_\mu[u]$ (where $g_\mu$ is taken in the sense of
distributions) is the unique entropy solution of (FDE).
\end{theorem}
\begin{proof}
For the HJB equation well-posedness of viscosity solutions for $\eps\geq0$ and the uniform convergence $u^\eps\ra u^0$ is fairly standard and can be found
e.g.~as a simple special case of results in \cite{JK:ContDep}.

Existence and uniqueness in (ii) follow from this paper for $\eps=0$. The arguments in this paper can easily be extended to include the $\eps\Delta v$-term
(this is standard) and hence we have (ii) for any $\eps>0$. The limit $v^\eps\ra v$ can be obtained through a standard Kuznetzov type argument,
cf.~\cite{Alibaud,Cifani/Jakobsen/Karlsen} for the case when $A$ is linear. We will give the result for the non-linear case in a future paper. The regularity
for $\eps>0$ is clear since the $\eps\Delta v$-term is the principal term in the equation. It follows e.g.~from (i) energy estimates and a classical parabolic
compactness argument (yields $L^2(0,T;H^1(\R^d))$-solutions) and (ii) regularity theory for the Heat equation combined with bootstrapping (yields smooth
solutions). The fractional term is always related to integer order derivatives through interpolation estimates. The detailed proof is long and rather classical
and is best left to the interested reader.
\end{proof}

\begin{remark}
Such correspondences between HJB equations and degenerate convection
diffusion equation can be useful for at least two reasons.
\begin{itemize}
\item[1)] They allow for integral representation formulas for the
  solutions  of the degenerate convection diffusion equations via
  representation formulas for the solutions of the HJB equations. See
   e.g. chapter 3.4 in \cite{Ev:Book} for the case of one dimensional scalar
   conservation laws.
\item[2)] They allow for efficient numerical methods for the
  non-divergence form HJB equation, by solving the divergence
  form degenerate convection diffusion equation by finite
  elements or spectral methods and then using the correspondence (and
  the HJB equation) to
  find the HJB solution.
\end{itemize}
The solutions of the above HJB equations are value functions of suitably defined stochastic differential games (see e.g. \cite{FS:Book}), i.e.~they have
integral representation formulas. Since HJB equations are fully non-linear non-divergence form equations, it is not natural or easy to solve them directly by
well-established, flexible, and efficient methods like the finite element and spectral methods. Such methods do apply to divergence form equations like the
degenerate convection diffusion equations (cf. e.g. \cite{CS98,JL08,Cifani/Jakobsen/Karlsen2}).
\end{remark}

\section{Numerical experiments}
We conclude this paper by presenting some experimental results
obtained using the numerical method \eqref{scheme} with $d=1$. We
simulate fractional strongly degenerate equations and compare them to
fractional conservation laws and local convection diffusion
equations. Our simulations give some insight into how the solutions of
these new equations behave. Note that this type of fractional equations
have never been simulated (or analyzed) before.

In our computations, we restrict ourselves to the bounded region
$\Omega=\{x:|x|\leq 2\}$ and impose zero Dirichlet boundary conditions
on the whole exterior domain $\{x:|x|>2\}$. We consider the degenerate
fractional convection-diffusion equations with Burgers type convection
($f(u)=u^2/2$),
\begin{align}\label{burgers}
\partial_tu+u\partial_xu=g[A(u)],
\end{align}
and fractional degenerate diffusion equations ($f\equiv0$),
\begin{align}\label{pure}
\partial_tu=g[A(u)],
\end{align}
for two different strongly degenerate diffusions, defined through
 two different $A$'s:
\begin{equation*}
\begin{split}
A_1(u)=\max(u,0)
\end{split}
\end{equation*}
and
\begin{equation*}
\begin{split}
A_2(u)&=\left\{
\begin{array}{ll}
0&u\leq0.5,\\
5(2.5u-1.25)(u-0.5)&0.5<u\leq0.6,\\
1.25+2.5(u-0.6)&u>0.6.
\end{array}
\right.
\end{split}
\end{equation*}

The numerical experiments below show e.g.~how solutions of \eqref{1} can develop shock discontinuities in finite time for all $\lambda\in(0,2)$. Furthermore,
they show that, contrary to the linear case, equation \eqref{pure} does not have smooth solutions for $t>0$ when the initial data is non-smooth. We also
observe that for $\lambda\approx2$, solutions are very close to solutions of the corresponding local problem with $\lambda=2$.

\begin{figure}[h]
\subfigure[]{
\includegraphics[width=60mm,height=40mm]{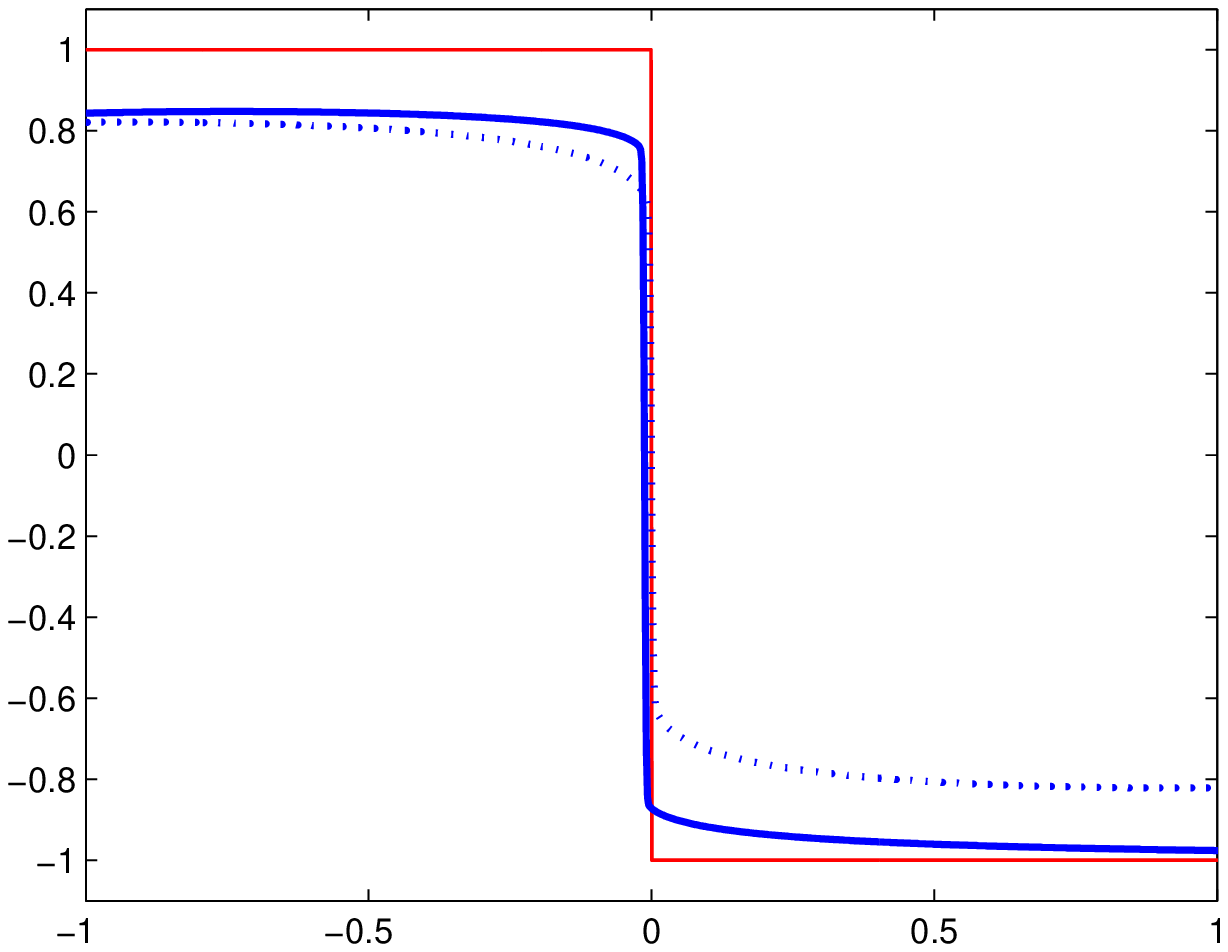}}
\subfigure[]{
\includegraphics[width=60mm,height=40mm]{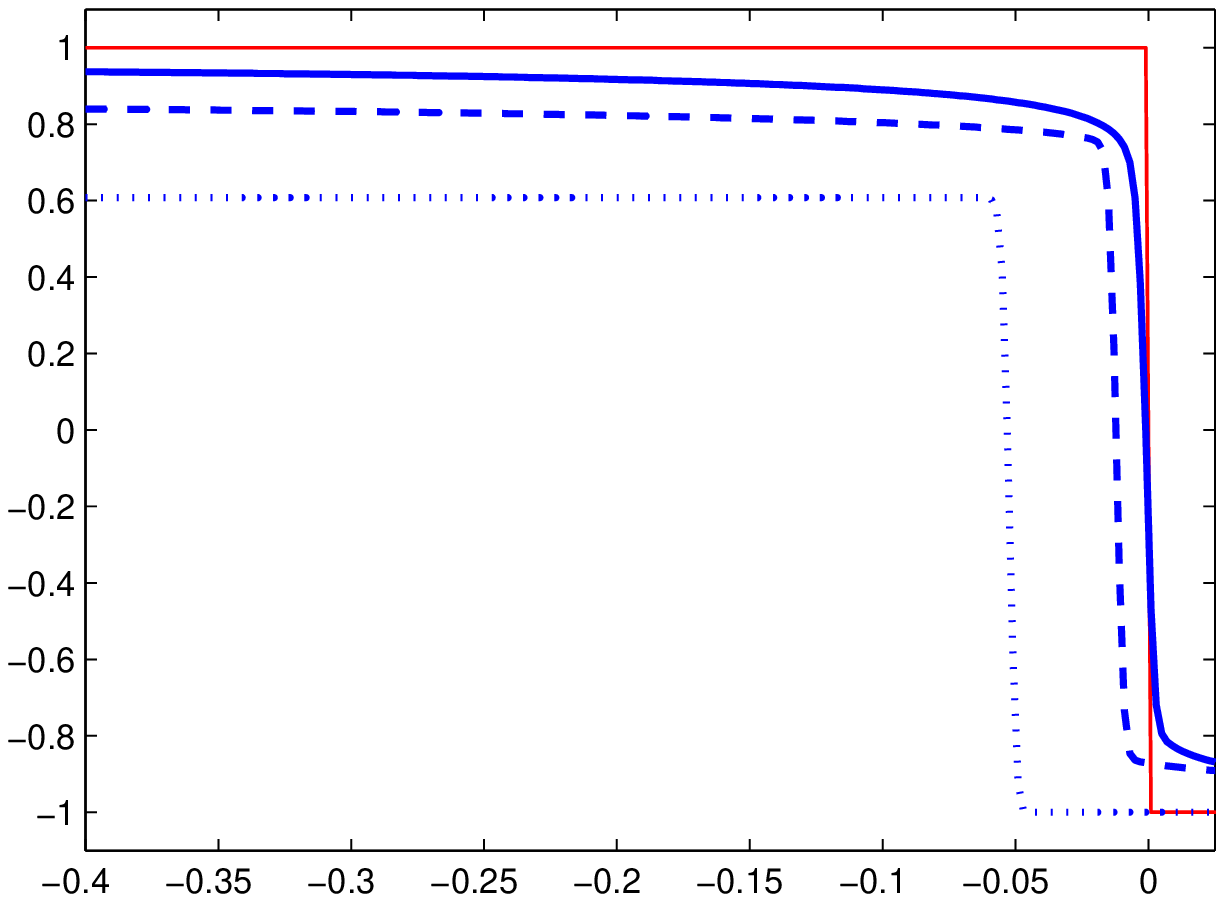}}
\caption{Numerical solutions of \eqref{burgers} at $T=0.5$ with $\Delta x=1/500$ and piecewise constant initial data: (a)  $A=A_1$
(solid) and $A(u)=u$ (dotted) for $\lambda=0.5$; (b) $A=A_1$ with $\lambda=0.001$ (dotted), $\lambda=0.5$ (dashed), and $\lambda=0.999$
(solid).}\label{fig_a}
\end{figure}

In figure
 \textsc{Figure} \ref{fig_a} (a)--(b) we plot the solutions of
 \eqref{burgers} with linear and non-linear fractional diffusion
 ($A(u)=u$
and $A=A_1$) to show how the non-linearity influences both the shock size and speed.

\begin{figure}[h]
\subfigure[]{
\includegraphics[width=60mm,height=40mm]{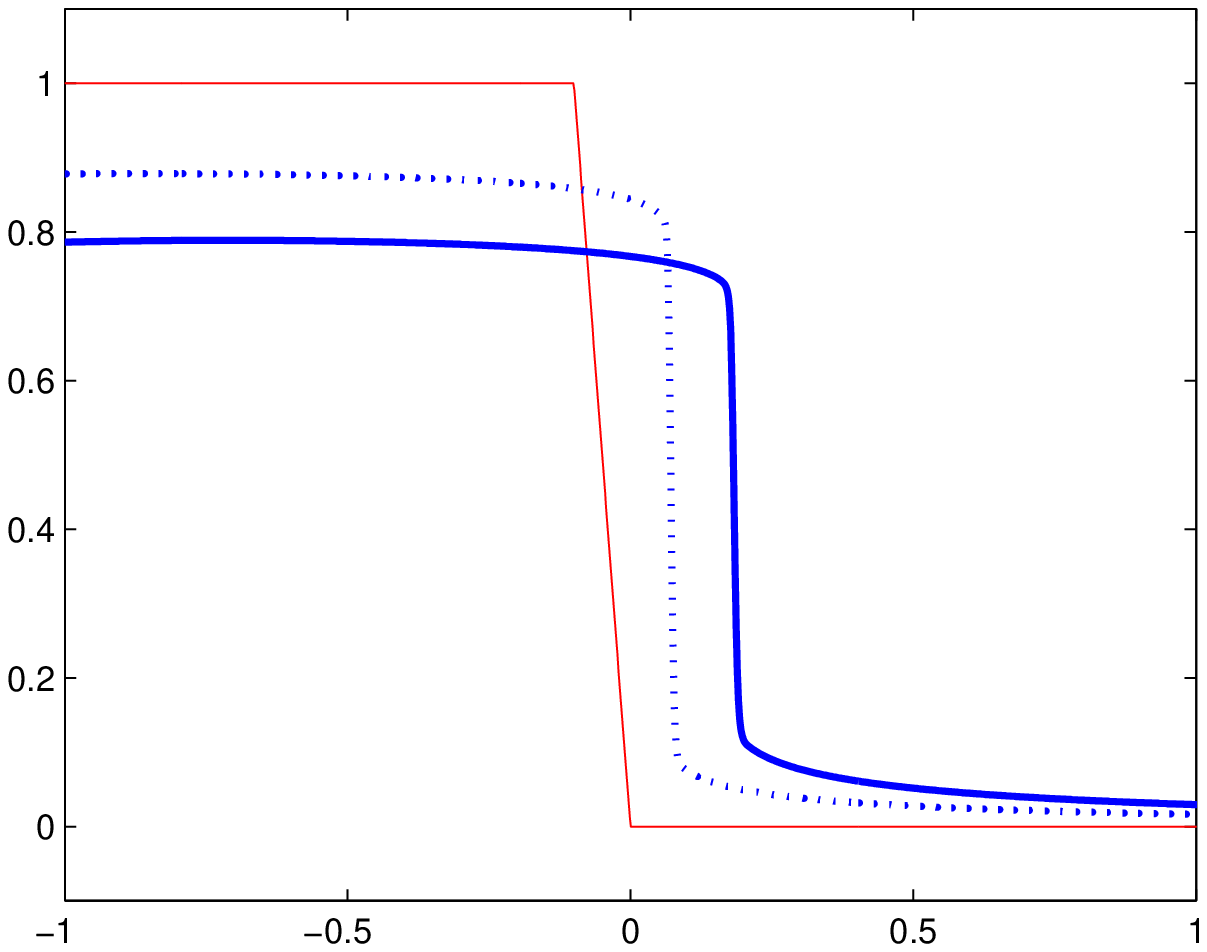}}
\subfigure[]{
\includegraphics[width=60mm,height=40mm]{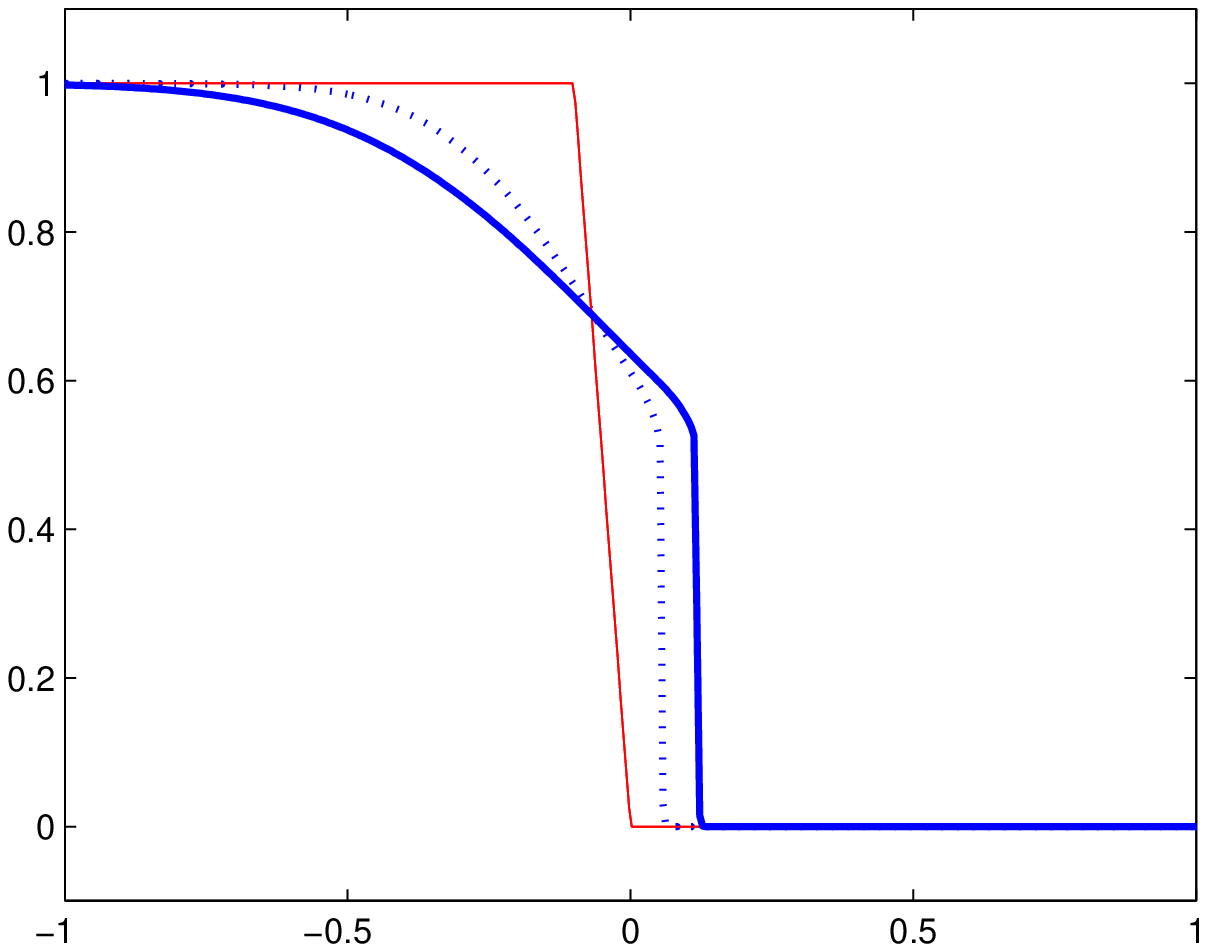}}
\caption{Burgers's flux and $A=A_2$ with $\Delta x=1/500$ and piecewise linear initial data: (a) solutions of \eqref{scheme} with $\lambda=0.3$ at $T=0.25$
(dotted) and $T=0.5$ (solid); (b) solutions of \eqref{scheme_conv} at $T=0.01$ (dotted) and $T=0.025$ (solid).}\label{fig_b}
\end{figure}

\textsc{Figure} \ref{fig_b} (a) shows that a shock discontinuity develops  in finite time in the region where $A_2$ is zero. This phenomenon is well known for
degenerate convection-diffusion equations \eqref{conv_diff} as shown
in \textsc{Figure} \ref{fig_b} (b). Here and in what follows, we have used the convergent
numerical scheme (cf.~\cite{Evje/Karlsen})
\begin{equation}\label{scheme_conv}
\begin{split}
U_{i}^{n+1}=U_{i}^{n}-\Delta tD^-F(U_{i}^{n},U_{i+1}^{n})+(2\pi)^2\Delta tD^-\left(\frac{A(U^n_{i+1})-A(U^n_i)}{\Delta x}\right).
\end{split}
\end{equation}
to compute the solutions of degenerate convection-diffusion equations \eqref{conv_diff}.

\begin{figure}[h]
\subfigure[]{
\includegraphics[width=60mm,height=40mm]{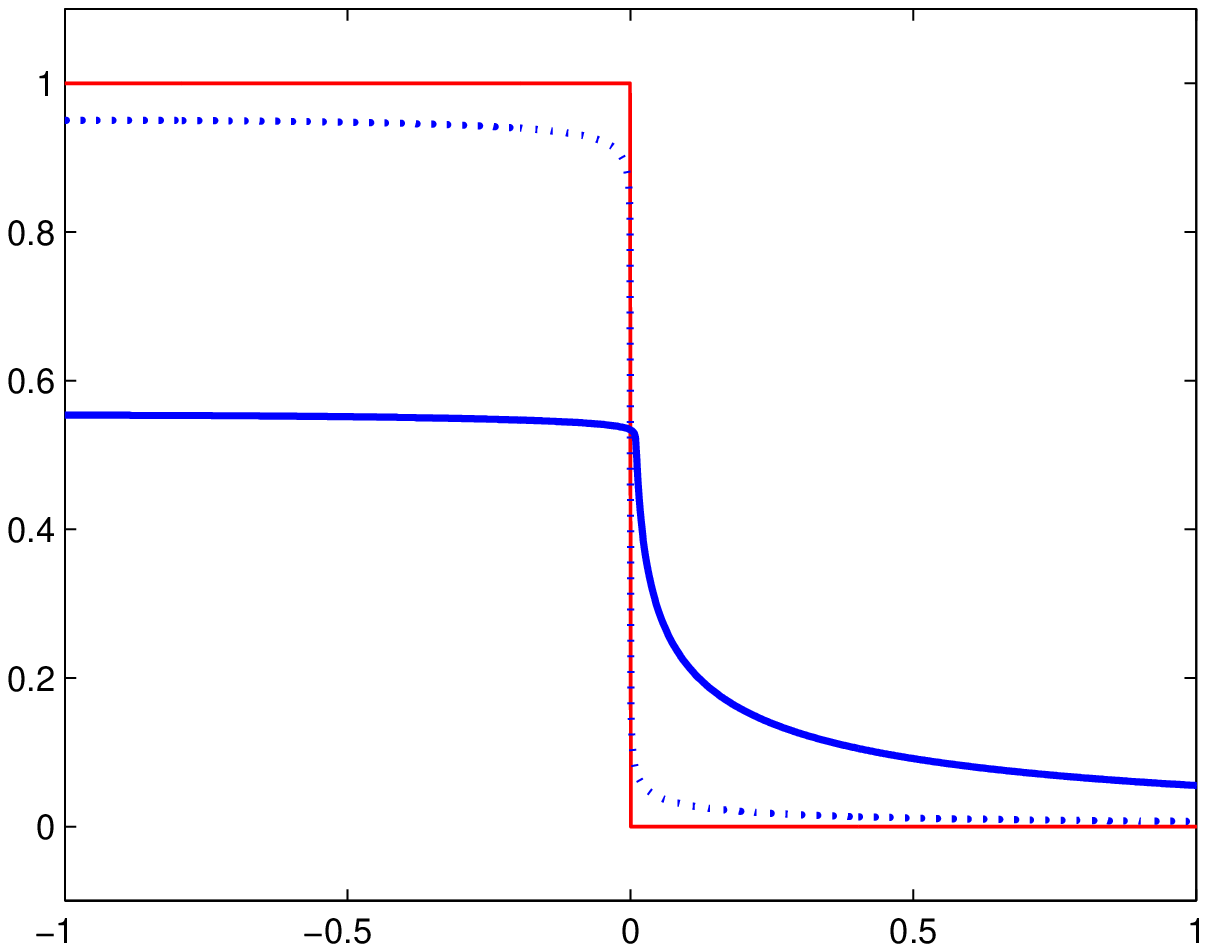}}
\subfigure[]{
\includegraphics[width=60mm,height=40mm]{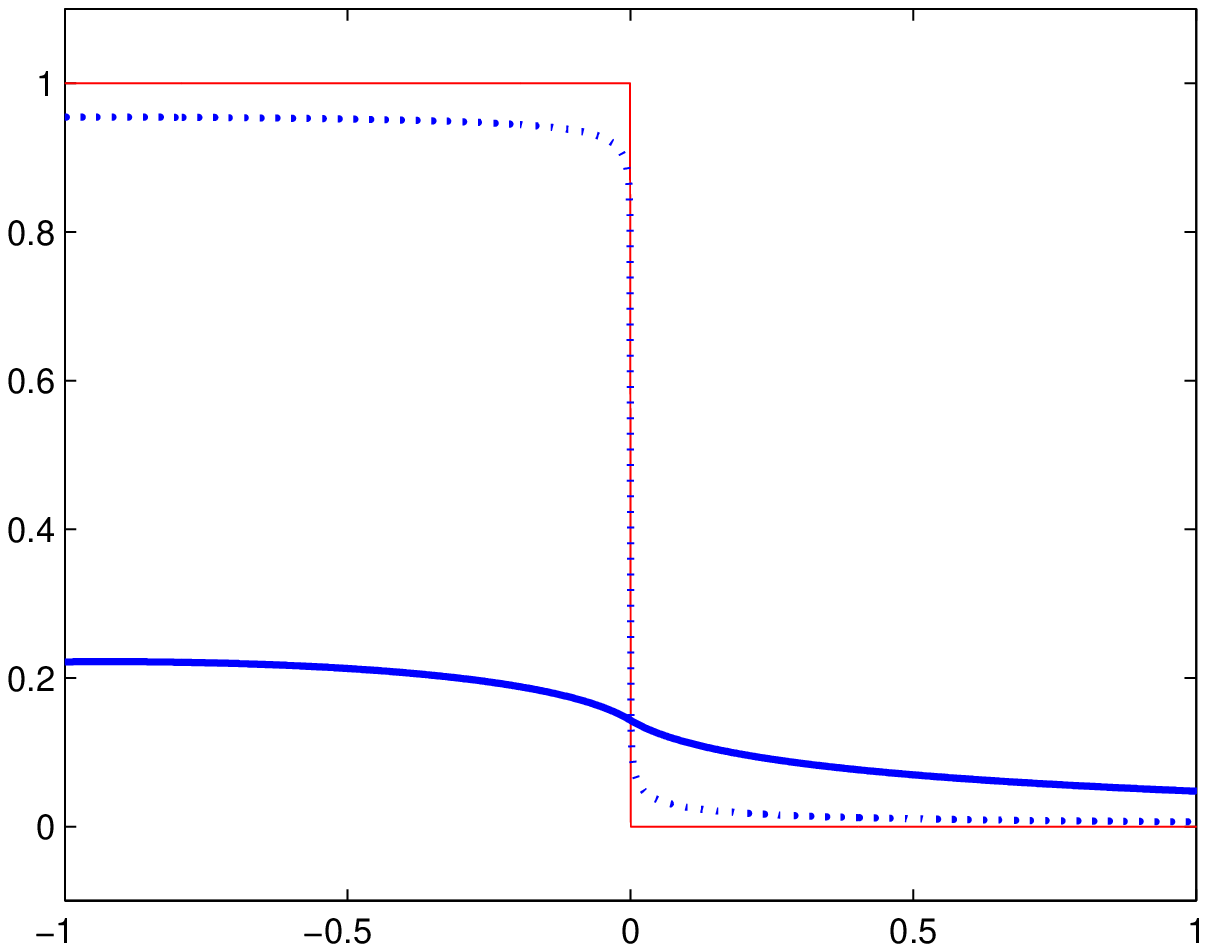}}
\caption{Solutions of \eqref{pure} with $\Delta x=1/500$, $\lambda=0.3$, and piecewise constant initial data: (a) $A=A_1$ with $T=0.1$ (dotted) and $T=3$
(solid); (b) $A(u)=u$ with $T=0.1$ (dotted) and $T=3$ (solid).}\label{fig_c}
\end{figure}

\textsc{Figure} \ref{fig_c} (a) displays the solutions of \eqref{pure} with $A(u)=u$ and $A=A_2$. Note that, when $A=A_2$, the initially discontinuous solution
becomes continuous in finite time but not differentiable. In the non-degenerate case, $\partial_t u=g[u]$, the initially discontinuous solution becomes smooth
immediately for all values of $\lambda$, cf.~\textsc{Figure} \ref{fig_c} (b). This behavior agrees with results from \cite{Droniou/Gallouet/Vovelle}.

\begin{figure}[t]
\subfigure{
\includegraphics[width=60mm,height=40mm]{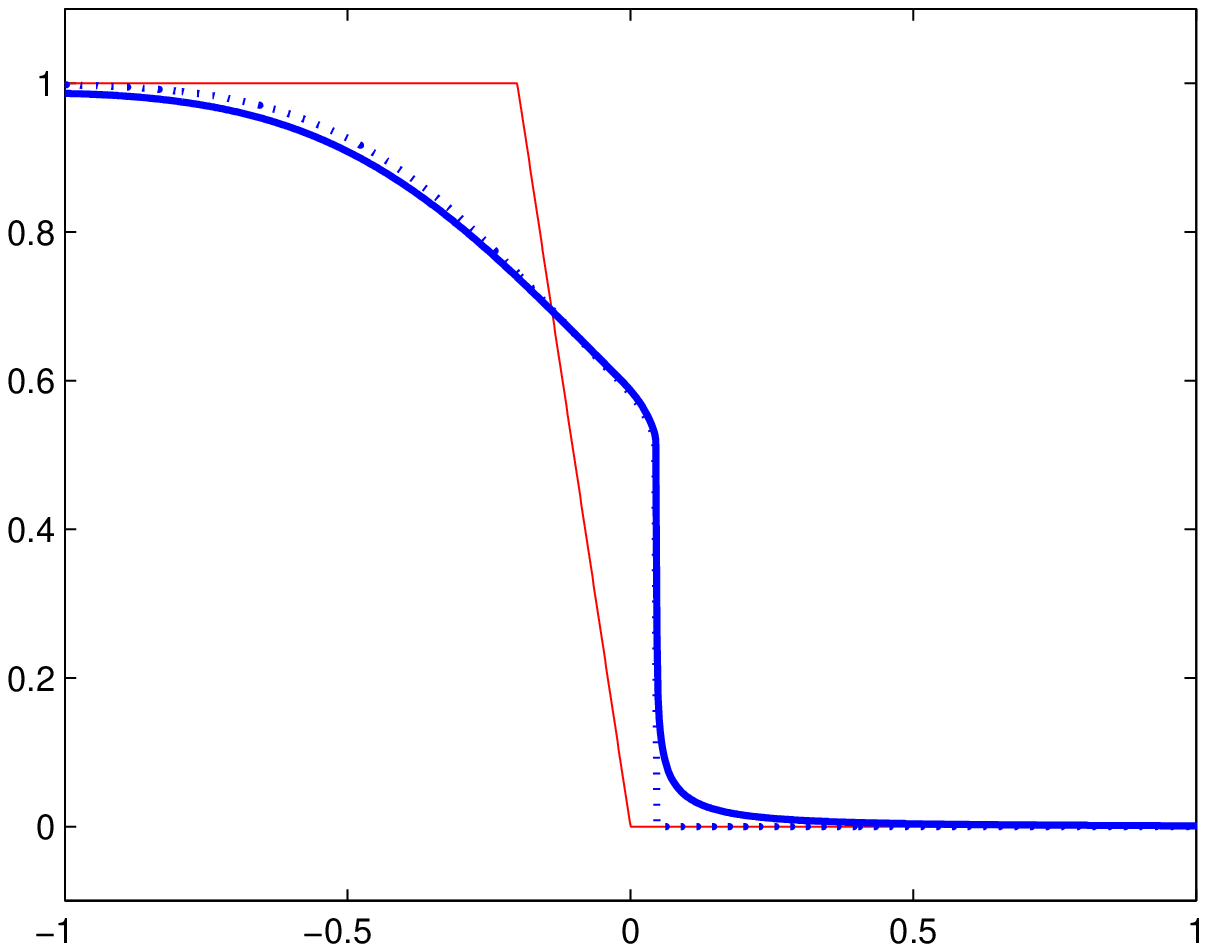}}
\subfigure{
\includegraphics[width=60mm,height=40mm]{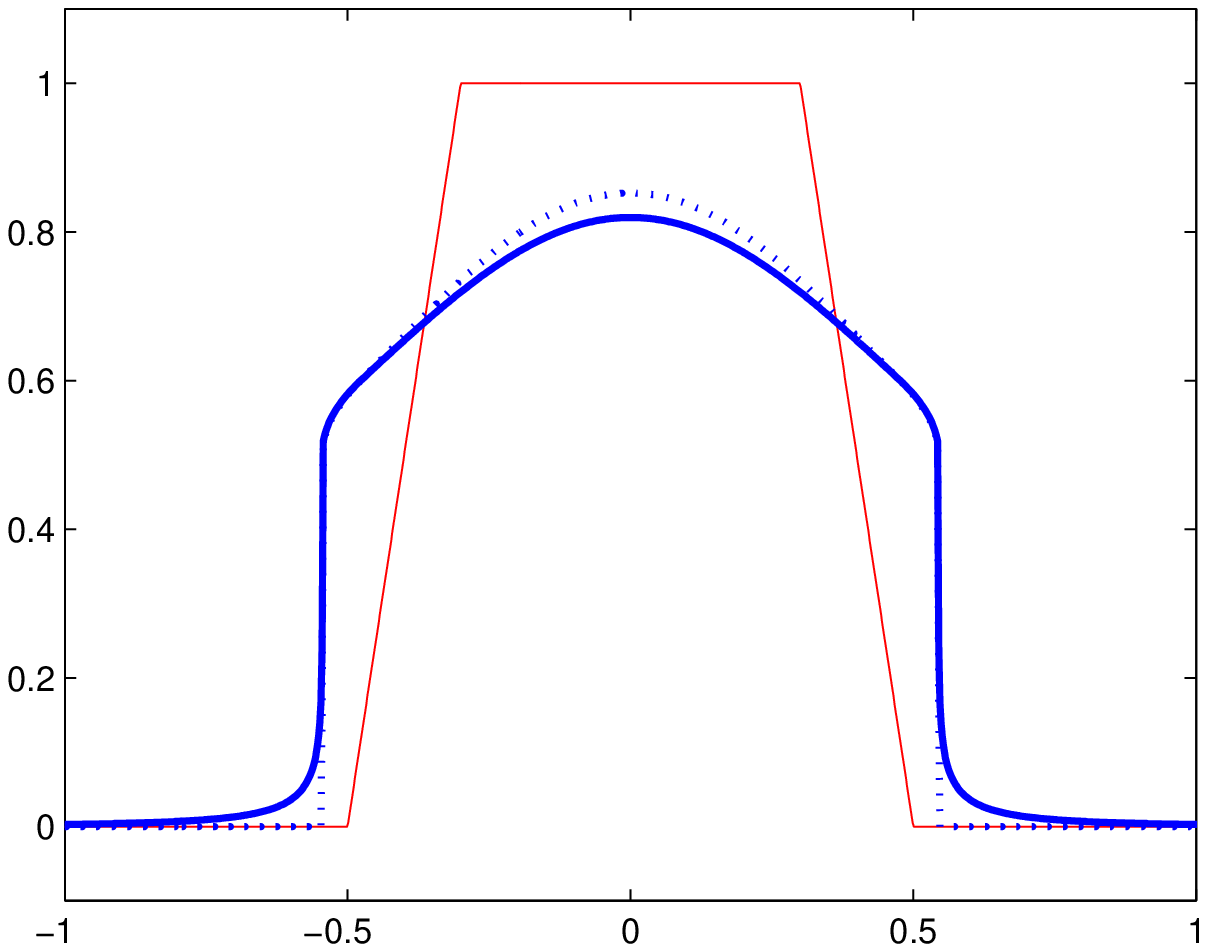}}
\caption{$A(u)=A_2$ with $T=0.005$, $\Delta x=1/500$, and piecewise constant initial data: (a)--(b) solutions of the non-local numerical method \eqref{scheme}
(solid) with $\lambda\approx2$ compared with solutions of the local numerical method \eqref{scheme_conv} (dotted).}\label{fig_AAA}
\end{figure}

In \textsc{Figure} \ref{fig_AAA} we compare the solutions of \eqref{pure} for $\lambda\approx2$, with the solutions of a properly scaled equation
\eqref{conv_diff} ($\lambda=2$). We use our scheme \eqref{scheme} to compute the first set of solutions, while scheme \eqref{scheme_conv} is used to compute
the second. Again, we have restricted our computational domain to $\Omega$. As expected, the solutions of the two equations are very close since
$-(-\Delta)^{\frac\lambda2}\phi\ra \Delta\phi$ as $\lambda\ra2$ for regular enough $\phi$. The two methods are however fundamentally different:
\eqref{scheme_conv} uses a three-points stencil, while \eqref{scheme} uses a whole-domain stencil.

%%%%%%%%%%%%%%%%%%%%%%%%%%%%%%%%%%%%%%%%%%%%%%%%%%%%%%%%%%%%%%%%%%%%%%%%%%%%%%%%%%%%%%%%%%%%%%%%%%%%%%%%%%%%%%%%%%%%%%%%%%%%%%%%%%%%%%%%%%%%%%%%%%%%%%%%%%%%%%%%%%%%%%%%%%%
%%%%%%%%%%%%%%%%%%%%%%%%%%%%%%%%%%%%%%%%%%%%%%%%%%%%%%%%%%%%%%%%%%%%%%%%%%%%%%%%%%%%%%%%%%%%%%%%%%%%%%%%%%%%%%%%%%%%%%%%%%%%%%%%%%%%%%%%%%%%%%%%%%%%%%%%%%%%%%%%%%%%%%%%%%%

\appendix
\section{A technical result}
In this section, we prove a technical result used in the proof of Lemma \ref{interp}.
\begin{lemma}\label{lem:BV}
Let $u\in BV(\R^d)$, then
\begin{equation}\label{kkk}
\begin{split}
\int_{\mathbb{R}^d}|u(x+z)-u(x)|\ \dif x\leq \sqrt{d}\, |z||u|_{BV(\mathbb{R}^d)}.
\end{split}
\end{equation}
\end{lemma}
Note that a more refined argument would give a factor $1$ instead of $\sqrt{d}$ in \eqref{kkk}. This is unimportant in this paper and we skip it. We now give a
proof for \eqref{kkk} in the case $d=2$, analogous ideas can then be used in higher dimensions.
\begin{proof}
We define the total variation $|u|_{BV(\mathbb{R}^2)}$ as, cf.~\cite[expression A.19]{Holden/Risebro},
\begin{equation}\label{def_multi}
\begin{split}
|u|_{BV(\mathbb{R}^2)}=\int_{\mathbb{R}}|u(x_1,\cdot)|_{BV(\mathbb{R})}\ \dif x_1+\int_{\mathbb{R}}|u(\cdot,x_2)|_{BV(\mathbb{R})}\ \dif x_2.
\end{split}
\end{equation}
Then, since $\int_{\R}|u(x+z)-u(x)|\ \dif x\leq |z||u|_{BV(\R)}$, we write
\begin{equation*}
\begin{split}
\int_{\mathbb{R}^2}|u(x+z)-u(x)|\ \dif x&=\int_{\mathbb{R}^2}|u(x_1+z_1,x_2+z_2)-u(x_1,x_2)|\ \dif x_1\dif x_2
\end{split}
\end{equation*}
which, by triangle inequality, is less than or equal to
\begin{equation*}
\begin{split}
&\int_{\mathbb{R}^2}|u(x_1+z_1,x_2+z_2)-u(x_1,x_2+z_2)|\ \dif x_1\dif x_2\\
&\qquad+\int_{\mathbb{R}^2}|u(x_1,x_2+z_2)-u(x_1,x_2)|\ \dif x_1\dif x_2\\
&\leq|z_1|\int_{\mathbb{R}}|u(\cdot,x_2+z_2)|_{BV(\mathbb{R})}\ \dif x_2\\
&\qquad+|z_2|\int_{\mathbb{R}}|u(x_1,\cdot)|_{BV(\mathbb{R})}\ \dif x_1\\
&\leq \sqrt{2}|z||u|_{BV(\mathbb{R}^2)},
\end{split}
\end{equation*}
thanks to \eqref{def_multi} and inequality $|z_1|+|z_2|\leq \sqrt{2}\,|z|$.
\end{proof}

\section*{Acknowledgement}
We would like to thank Natha\"{e}l Alibaud, Harald Hanche-Olsen, and Boris Andreianov for many helpful discussions concerning this paper. We would also like to
thank the two anonymous referees for their very careful reports. All these people have helped us improve this paper a lot.

\end{document}